\documentclass[10pt]{article}
\usepackage{amsmath, amssymb,amsxtra, graphicx, float}
 \pdfoutput=1
\newtheorem{thm}{{\sc Theorem}}[section]

\newtheorem{df}[thm]{{\sc Definition}}

\newcommand{\LR}{\hbox{Little\-wood-Richard\-son}}

\catcode`\@=11 \font\linear=line10 scaled \magstep5
\def\slant{{\linear ,}}
\def\young@lign{\everycr{}\tabskip0pt\halign}
\def\Mathstrut@{\setbox\z@\hbox{$($}\setbox\tw@\null\ht\tw@\ht\z@\dp\tw@\dp\z@
 \box\tw@}
\def\b@m#1{$\m@th\underline{#1}$}
\def\t@p#1{$\m@th\overline{#1}$}

\def\young{
\setbox\strutbox=\hbox{\vrule height3pt depth3.5pt width\z@}
 \offinterlineskip%
 {}\,\vcenter
\bgroup
     \def\1{{}}
     \def\2{\1&\1}
     \def\3{\1&\1&\1}
     \def\4{\1&\1&\1&\1}
\let\\=\cr
 \tabskip0pt\baselineskip0pt\m@th
 \young@lign
 \bgroup\vrule\b@m{\hbox to .75 em{\strut\hfil$##$\hfil}}\vrule %
   &&\b@m{\hbox to .75em{\strut\hfil$##$\hfil}}\vrule\crcr
   \noalign{\hrule}
  }

\def\endyoung{\egroup \egroup\,}
\def\slyoung{
   \setbox\strutbox=\hbox{\vrule height10pt depth2pt width\z@}
   \offinterlineskip%
   {}\,\vcenter
\bgroup
     \def\1{{}}
     \def\2{\1&\1}
     \def\3{\1&\1&\1}
     \def\4{\1&\1&\1&\1}
\let\\=\cr
 \tabskip0pt\baselineskip0pt\m@th
 \young@lign
 \bgroup\lower2pt\hbox{\slant}\kern-25pt\b@m{\hbox to 2.5em{\strut\hfil$##\;\;$}}\vrule
   &&\b@m{\hbox to 2.5em{\strut\hfil$##$\hfil}}\vrule\crcr
   \noalign{\hrule}
  }

\def\endslyoung{\egroup \egroup\,}
\def\frame #1#2#3#4{\vbox{\hrule height #1pt
 \hbox{\vrule width #1pt\kern #2pt
 \vbox{\kern #2pt
 \vbox{\hsize #3\noindent #4}
 \kern #2pt}
 \kern #2pt\vrule width #1pt}
 \hrule height0pt depth #1pt}}

\def\nframe #1#2#3#4{\vbox{
 \hrule height #1pt width0pt
 \hbox{\vrule height0pt width #1pt\kern #2pt
 \vbox{\kern #2pt
 \vbox{\hsize #3\noindent #4}
 \kern #2pt}
 \kern #2pt\vrule width #1pt height0pt}
 \hrule height0pt width0pt}}

\def\leftwall #1#2#3#4{\vbox{
 \hrule height #1pt width0pt
 \hbox{\vrule width #1pt\kern #2pt
 \vbox{\kern #2pt
 \vbox{\hsize #3\noindent #4}
 \kern #2pt}
 \kern #2pt\vrule width #1pt height0pt}
 \hrule height0pt width0pt}}

\def\rightwall #1#2#3#4{\vbox{
 \hrule height #1pt width0pt
 \hbox{\vrule height0pt width #1pt\kern #2pt
 \vbox{\kern #2pt
 \vbox{\hsize #3\noindent #4}
 \kern #2pt}
 \kern #2pt\vrule width #1pt }
 \hrule height0pt width0pt}}

\def\b #1{\frame{.3}{2}{8pt}{\centerline{#1}\vphantom{(}}}

\def\e{\vphantom{e}} 
\def\f{\b{\e}} 
\def\nl{\hfill\break} 

\usepackage{amsmath}
\usepackage{float}

\begin{document}

\noindent {\bf Flows on Honeycombs and Sums of Littlewood-Richardson Tableaux}

\medskip

\noindent
Glenn D. Appleby\\
 Tamsen Whitehead\\
{\em Department of Mathematics\\
and Computer Science,\\
Santa Clara University\\
Santa Clara,  CA 95053}\\
gappleby@scu.edu, tmcginley@scu.edu

\noindent \mbox{} \hrulefill \mbox{}
\begin{abstract}
 Suppose $\mu$ and $\mu'$ are two partitions. We will let $\mu \oplus \mu'$ denote the {\em direct sum} of the partitions, defined as the sorted partition made of the parts of $\mu$ and $\mu'$. In this paper, we define a summation operation on two \LR\ fillings of type $(\mu, \nu;\lambda)$ and $(\mu', \nu';\lambda')$, which results in a \LR\ filling of type $(\mu\oplus \mu', \nu\oplus \nu' ;\lambda\oplus \lambda')$.  We give an algorithm to produce the sum, and show that it terminates in a \LR\ filling by defining a bijection between a \LR\ filling and a {\em flow} on a honeycomb, and then showing that the overlay of the two honeycombs of appropriate type corresponds to the sum of the two fillings.
\end{abstract}

\noindent \mbox{} \hrulefill \mbox{}

\section{Introduction}
There has been an active interest in relating the combinatorics of \LR\ fillings to other mathematical objects. Survey papers by Fulton \cite{fulton} and Zelevinsky~\cite{zel} demonstrate that these combinatorial objects appear in a wide variety of contexts including representation theory, the eigenvalue structure of Hermitian matrices, and the Schubert calculus. In particular, the relationship between \LR\ fillings and combinatorial invariants called {\em honeycombs} and {\em hives} has been investigated (definitions and examples provided below).

In this paper we show how to construct a canonically defined {\em flow} on a given honeycomb, which is decomposed into parts that also determine the parts of the \LR\ filling corresponding to the honeycomb. We then construct an algorithm on a pair of \LR\ fillings that produces a ``sum'' of the filling. The overlay the two honeycombs corresponds to the sum of two fillings and, it is shown, the overlay of the associated canonical flows, which is typically not canonical. We demonstrate that our algorithm on \LR\ fillings corresponds to ``resolving'' the overlaid flow into its canonical flow, where swapping parts of the decomposition of the flow parallels the steps of our algorithm on the parts of the filling. Sums of honeycombs and/or \LR\ fillings are connected to open questions involving, among other things, spectra of sums of Hermitian matrices, and our algorithm is part of a larger program to make these connections more explicit.

Knutson and Tao's original work~\cite{knut} on honeycombs provided a framework to relate \LR\ fillings to problems concerning the spectra of sums of Hermitian matrices. Their work helped complete the classification of those triples of sequences of real numbers $(\mu, \nu, \lambda)$ for which there exist $k \times k$  Hermitian matrices $M$ and $N$ such that the spectrum of $M$ is $\mu$, the spectrum of $N$ is $\nu$, and the spectrum of $M+N$ is $\lambda$. Knutson and Tao pointed out that the overlay of honeycombs {\em should} correspond to the direct sum of Hermitian matrices. That is, there {\em should} be a way to determine a honeycomb ${\cal H}$ of type $(\mu, \nu; \lambda)$ associated to a $k \times k$ Hermitian pair $M$ and $N$, and another honeycomb ${\cal H}'$ associated to a $\ell \times \ell$ pair $M'$ and $N'$, such that the overlay of ${\cal H}$ and ${\cal H'}$ is the honeycomb associated to the $(k + \ell) \times (k+\ell)$ pair $M \oplus M'$ and $N \oplus N'$. While it is possible to verify this in very simple cases (such as diagonal matrices) this goal has not been realized generally.

Subsequent work~\cite{KTT,knut-purb} highlighted the overlay construction of honeycombs in order to answer questions regarding factorization of \LR\ coefficients. Our paper, however, addresses some of the original questions for which honeycombs were applied in the context of matrix decompositions. We provide here a combinatorial algorithm which directly constructs a \LR\ filling as a sum of two other such fillings, and we prove that this algorithm precisely matches the overlay construction of the associated honeycombs. The algorithm is of independent combinatorial interest, and it also helps us understand matrix decompositions in the Hermitian, and other, contexts.

Let $c_{\mu,\nu}^{\lambda}$ denote number of the classical \LR\ fillings with skew-shape $\lambda / \mu$ with content $\nu$. In this paper, a {\em partition} will denote a non-increasing sequence of numbers (typically integers). If $\mu = (\mu_{1}, \ldots , \mu_n)$, then some $\mu_i$ may be zero, and we will call $n$ (the number of terms in the sequence) the {\em length} of the partition. Given partitions $\mu^{(1)}$ and $\mu^{(2)}$ (not necessarily the same length), let $\mu^{(1)} \oplus \mu^{(2)}$ denote the partition obtained by sorting the parts of $\mu^{(1)}$ and $\mu^{(2)}$ together.  We will call this the {\em direct sum} of the partitions $\mu^{(1)}$ and $\mu^{(2)}$.  Suppose $\mu$, $\nu$ and $\lambda$ are partitions all of length $n$, so that $\mu = (\mu_1, \mu_2, \ldots , \mu_n)$, etc. Let $I = (i_1, i_2, \ldots , i_k)$ be a sequence of indices from $\{ 1, 2, \ldots n \}$, and similarly for length $k$ sets of indices $J=(j_{\ell})$ and $K=(\kappa_{\ell})$. Let $\mu_{I} = (\mu_{i_{1}}, \mu_{i_{2}}, \ldots , \mu_{i_{k}} )$ be the associated length $k$ sub-sequence of $\mu$, and let $\mu_{I'}$ denote the complementary sequence of parts of $\mu$ not appearing in $\mu_{I}$, and similarly for $\nu_{J}, \nu_{J'}$ and $\lambda_{K}, \lambda_{K'}$. We also let
\[ | \mu_{I} | = \mu_{i_{1}} + \mu_{i_{2}} + \cdots + \mu_{i_{k}} . \]
Note that, with these definitions, $\mu = \mu_{I} \oplus \mu_{I'}$, etc.

In~\cite{knut} Knutson and Tao completed a proof that a triple of eigenvalues $(\mu, \nu; \lambda)$ form the spectra for Hermitian matrices $M$, $N$, and $L$ such that $M+N=L$ if and only if a collection inequalities of the form
 \[ | \lambda_{K}| \leq |\mu_{I}|+ |\nu_{J}| \]
 are true, along with the trace condition $|\mu|+|\nu| = | \lambda|$. The collection of index sets $(I,J,K)$ determining the inequalities we shall call {\em Horn triples}, and the inequalities they define we call {\em Horn inequalities}, since it was Horn~\cite{Horn} who first conjectured the above result, though these sequences of index sets go by many names in the literature. It was proved subsequently, however, that the Horn inequalities are not minimal, and the Hermitian matrix existence problem is implied by inequalities determined by subset of the Horn triples, now known as {\em essential triples} determining the {\em essential Horn inequalities}.

King, Tollu, and Toumazet~\cite{KTT} considered Horn and essential triples in order to study hives and puzzles (combinatorial invariants related to honeycombs). Their results implied that if a triple $(\mu, \nu, \lambda)$ satisfied all the Horn inequalities, and if a Horn triple $(I,J,K)$ could be found so that the resulting inequality was tight:

\[ | \lambda_{K}| = |\mu_{I}|+ |\nu_{J}|, \]
then any honeycomb of type $(\mu, \nu, \lambda)$ was an overlay of a honeycomb of type $(\mu_{I}, \nu_{J}, \lambda_{K})$, and one of type $(\mu_{I'}, \nu_{J'}, \lambda_{K'})$.  If the Horn triple $(I,J,K)$ was {\em essential}, King, Tollu, and Toumazet proved the stronger statement that the collection of all honeycombs of type $(\mu, \nu, \lambda)$ decomposed in such a way as to imply the \LR\ coefficient
$c_{\mu, \nu}^{\lambda}$ factors as
$$c_{\mu \nu}^{\lambda} = c_{\mu_{I} \oplus \mu_{I'} ,\nu_{J} \oplus \nu_{J'}}^{\lambda_{K} \oplus \lambda_{K'}}=
c_{\mu_{I} \nu_{J}}^{\lambda_{K}}\cdot c_{\mu_{I'} \nu_{J'}}^{\lambda_{K'}}.$$
We note that in order for a triple $(I,J,K)$ to be Horn triple, it is necessary, but not sufficient that
\begin{equation} i_{i} + i_{2} + \cdots + i_{k} + j_1 + j_2 + \cdots + j_k = \kappa_1 + \kappa_2 + \cdots + \kappa_k + k. \label{top-eq}\end{equation}

However, not all honeycomb decompositions arise from products of \LR\ coefficients, or even from Horn triples, and these other decompositions are of interest.

As stated above, it was conjectured in~\cite{knut} that the direct sum of Hermitian matrices corresponds to an overlay of honeycombs, though an explicit map from Hermitian matrices to honeycombs, or conversely from \LR\ fillings to Hermitian matrices (the {\em Hermitian Realization Problem}) has not been found. In recent results by the authors~\cite{me,lrcon,lrreal} that relate \LR\ coefficients to the invariant factors of products of matrices over valuation rings, a direct correspondence {\em has} been achieved, and a very close connection between direct sums of matrices and overlays of honeycombs exists.  However, the honeycomb decompositions appearing in general matrix decompositions are not, typically, of the sort analyzed by the results of King, Tollu, and Toumazet~\cite{KTT}.  For example, the honeycomb:

 \begin{figure}[H]
\centering{\includegraphics[scale=.3]{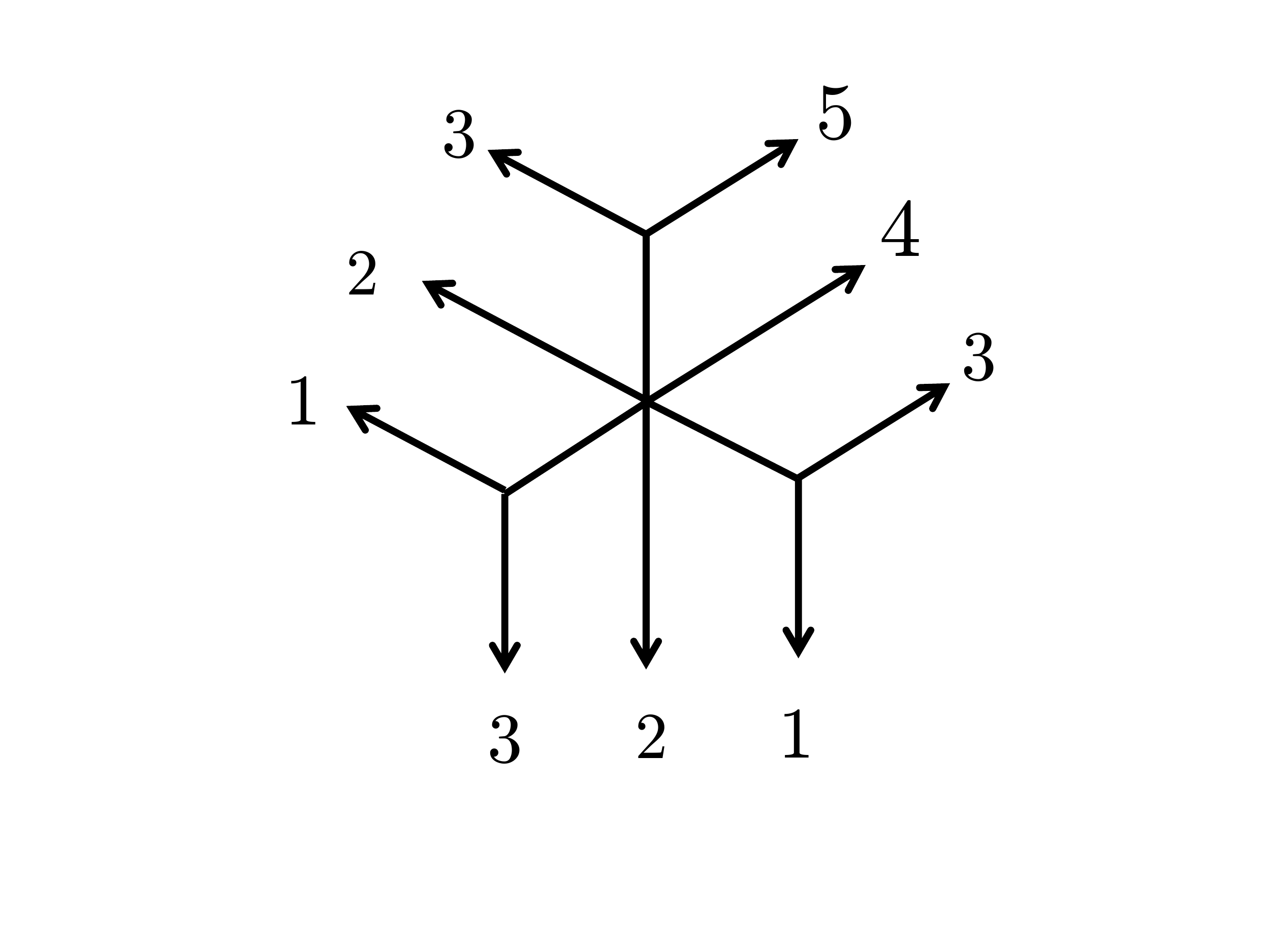}}
\caption{Honeycomb Example}
\label{Fig10}
\end{figure}

has a decomposition as an overlay in {\em two} inequivalent ways, first as a sum of two honeycombs:

\begin{figure}[H]
\centering{\includegraphics[scale=.45]{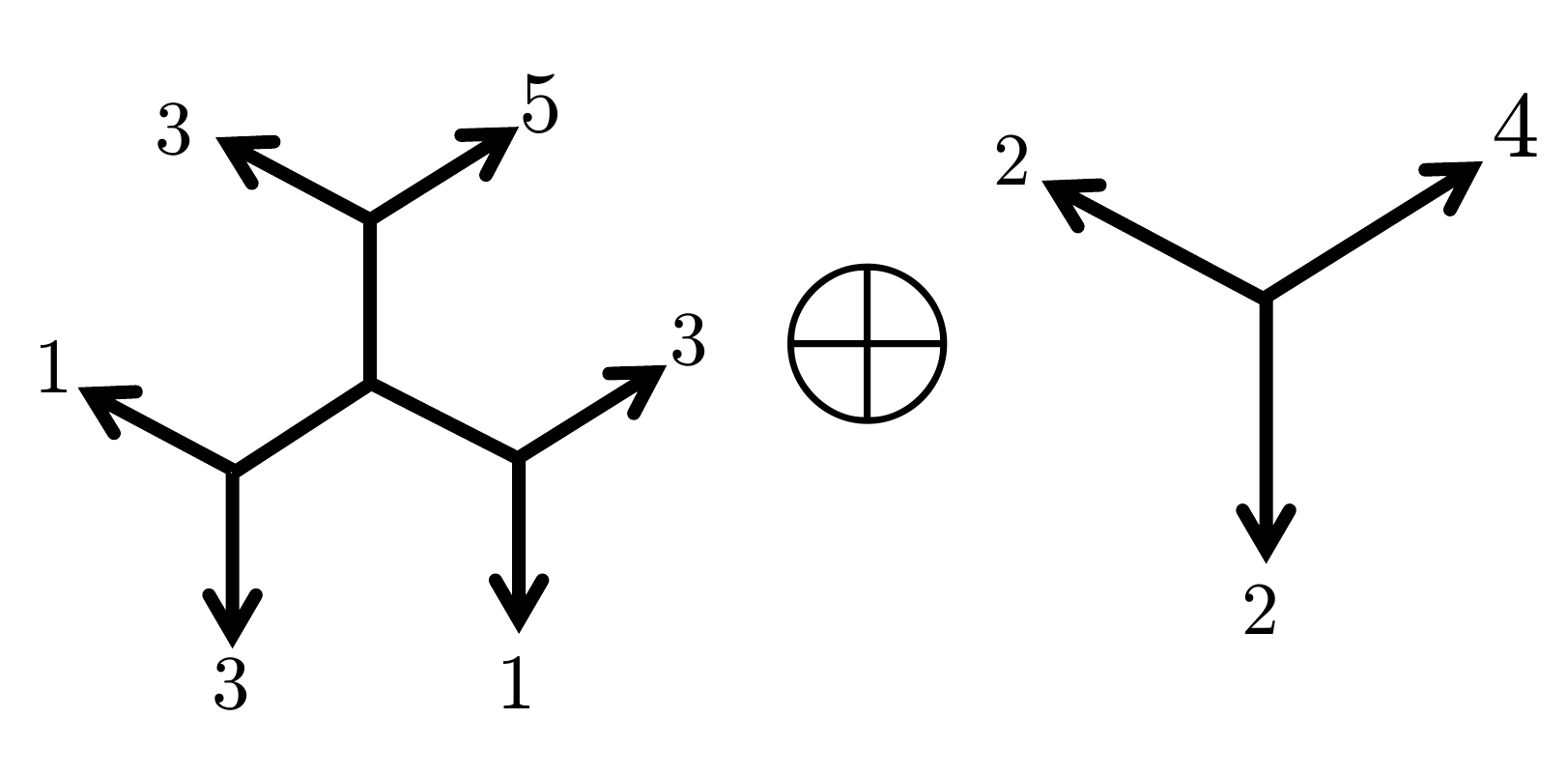}}
\caption{First Decomposition}
\label{Fig10}
\end{figure}
but also as a sum of three honeycombs:
\begin{figure}[H]
\centering{\includegraphics[scale=.45]{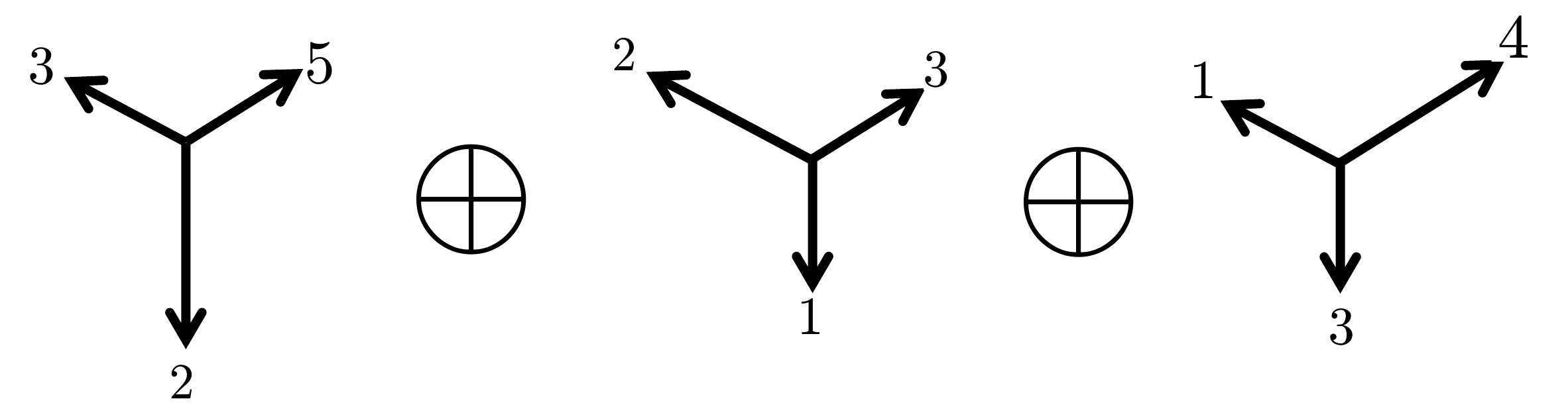}}
\caption{Second Decomposition}
\label{Fig10}
\end{figure}

 In neither decomposition do the index sets satisfy Equation \ref{top-eq} necessary for Horn or essential triples. The first decomposition uses the length-one index sets $((1),(2),(1)) \oplus ((2),(3),(3))\oplus ((3),(1),(2))$, while the second decomposition uses $( (3,1),(3,1),(3,1)) \oplus ((2),(2),(2))$.  What the above example shows is that the relation between \LR\ fillings and matrix (or honeycomb) decompositions is rather subtle, as there may be several inequivalent ways to realize a given honeycomb, with none determined by the Horn triples. The first decomposition may be built using diagonal matrices, but the second would not.

 While the overlay of two honeycombs is easy to construct, the resulting \LR\ filling formed by the overlay may be rather complicated. Conversely, it can be a challenge to detect whether a given honeycomb is an overlay of two others, and to determine the \LR\ coefficients of the two associated honeycombs can be even more difficult (and, as the above example shows, not uniquely determined).

Consequently, it is natural to ask if one can find a {\em sum} of \LR\ fillings of types
$(\mu, \nu; \lambda)$ and $(\mu', \nu'; \lambda')$ that will
result in a filling of type $( \mu\oplus \mu', \nu\oplus \nu';{\lambda\oplus \lambda'})$.
Below, we define such a sum algorithmically and show that it terminates in a \LR\ filling of the proper type.

 Also new is our use of {\em flows} on honeycombs. We define a bijection between \LR\ fillings of type $(\mu, \nu;\lambda)$ and flows on honeycombs of the same type. Our flow construction allows one to view the constraints appearing in the definition of \LR\ fillings as requirements on the crossing of flows along honeycomb paths. Using flows on hives (or their dual graphs) to determine the positivity of \LR\ coefficients has been useful for some problems~\cite{BH}, but our construction of flows on honeycombs is new and it greatly simplifies our analysis. It appears to provide a link between the combinatorial constructions of \LR\ coefficients and the variational definitions of eigenvalues and invariant factors.

\section{Notation and Definitions}


For two partitions $\alpha$ and $\beta$, we write $\beta \subseteq \alpha$, to
mean $\beta_{k} \leq \alpha_{k}$ for all $k \geq 1$.  This notation is
suggested by the fact that if we represent the partitions by decreasing,
left-justified rows of boxes (called the {\em diagram} or {\em Ferrers diagram}
of the partition), then $\beta \subseteq \alpha$ implies the diagram for
$\beta$ fits inside the diagram of $\alpha$. When $\beta \subseteq \alpha$, we
will denote by $\alpha / \beta$ the {\em skew shape} consisting of the
diagram of $\alpha$, with the diagram of $\beta$ removed.  In the example
below, $\alpha = (11, 10, 7,5)$,  $\beta =(7, 4, 2, 1)$ and and the skew shape $\alpha / \beta$
consists of the boxes of $\alpha$ containing an $x$.

\bigskip
\hspace{1.8in} \vbox{ \offinterlineskip \openup-1.5pt \nl
\f\f\f\f\f\f\f\b{$x$}\b{$x$}\b{$x$}\b{$x$} \nl
\f\f\f\f\b{$x$}\b{$x$}\b{$x$}\b{$x$}\b{$x$}\b{$x$}\nl
\f\f\b{$x$}\b{$x$}\b{$x$}\b{$x$}\b{$x$}\nl
\f\b{$x$}\b{$x$}\b{$x$}\b{$x$}\nl}
\bigskip

Our central combinatorial definition is the following:

\begin{df}
Let $\mu, \nu$, and $\lambda$ be partitions, with $length(\mu)$, $length(\nu) \leq
length(\lambda) \leq r$.  Let $S= (\lambda^{(0)},\lambda^{(1)}, \cdots ,
\lambda^{(r)})$ be a sequence of partitions in which $\lambda^{(0)} = \mu$, and $\lambda^{(i)} \subseteq \lambda^{(i+1)}$.
The sequence $S$ is called a {\em \LR\ Sequence} of type $(\mu, \nu;
\lambda)$ if there is a triangular array of non-negative integers $F=\{ k_{ij}: 1
\leq j \leq r,  j \leq i \leq r \}$ (called the {\em filling}) such
that, for $1 \leq i \leq j, \ \ \lambda^{(i)}_{j} = \mu_{j} + k_{1j}
+ \dots + k_{ij}$, subject to the conditions (LR1), (LR2), (LR3)
below. We shall say, equivalently, that any such set $F$
determines a {\em \LR\ filling} of the skew shape $\lambda / \mu$
with content $\nu$. \label{LR-def}

(LR1) (Sums) For all $1 \leq j \leq i \leq r$,
\[  \mu_{j}+\sum_{s=1}^{j} k_{sj} = \lambda_{j}, \quad \hbox{and} \quad \sum_{s=i}^{r}k_{is} =
\nu_{i}. \]

(LR2) (Column Strictness) For each $j$, for $2 \leq j \leq r$ and $1 \leq i \leq j$ we require
$\lambda^{(i)}_{j} \leq \lambda^{(i-1)}_{(j-1)},$ that is,
\[  \mu_{j}+ k_{1j} + \cdots k_{ij} \leq \,  \mu_{(j-1)} +k_{1,(j-1)} +
 \cdots + k_{(i-1),(j-1)} . \label{cs}
\]

(LR3)  (Word Condition) For all $1 \leq i \leq r-1$, $i \leq j \leq r-1$,
\[   \sum_{s=i+1}^{j+1} k_{(i+1),s} \  \leq \  \sum_{s=i}^{j}k_{is}. \] \label{wd}

Let $LR(\mu, \nu;\lambda)$ denote the set of \LR\ fillings of type $\lambda/\mu$ and content $\nu$.
\end{df}
This definition can be given a more visual interpretation.  If we take the skew diagram
$\lambda/\mu$, we can fill it with $\nu_1$ $1$'s, $\nu_2$ $2$'s, etc.
The first equality of (LR1) ensures that the sum of the number of boxes in row
$j$ of the filled diagram (including the empty boxes of the parts of $\mu$) is
$\lambda_{j}$, the $j$-th part of the partition $\lambda$, while in the
second equality we require that the sum of the number of $i$'s in all the rows
is $\nu_{i}$, the $i$-th part of the partition $\nu$.   (LR2) says that
the numbers in the filling are strictly increasing down columns. Lastly, (LR3)
indicates that the number of $i$'s in rows $i$ through $j$ is greater than or
equal to the number of $(i+1)$'s in rows $(i+1)$ through $(j+1)$.

For example, the diagram below is a \LR\ filling of $\lambda/\mu$ with content $\nu$,
where $\lambda = (11, 10, 7,5)$, $\mu = (7, 4, 2, 1)$ and $\nu= (8, 5, 4, 2)$. We call such a filled diagram a {\em \LR\ Tableau}.

\bigskip
\hspace{1.8in} \vbox{ \offinterlineskip \openup-1.5pt \nl
\f\f\f\f\f\f\f\b{$1$}\b{$1$}\b{$1$}\b{$1$} \nl
\f\f\f\f\b{$1$}\b{$1$}\b{$2$}\b{$2$}\b{$2$}\b{$2$}\nl
\f\f\b{$1$}\b{$2$}\b{$3$}\b{$3$}\b{$3$}\nl
\f\b{$1$}\b{$3$}\b{$4$}\b{$4$}\nl}

\bigskip
A consequence of (LR3) is that the filling cannot have any entry bigger than $i$ in row $i$.  So row 1 can only contain $1$'s, row 2 can only contain $1$'s and $2$'s and so on.

\begin{df} Given partitions $\mu$, $\nu$, and $\lambda$, we shall
let $c_{\mu \nu}^{\lambda}$ denote the number of \LR\ fillings of the skew shape $\lambda / \mu$ with content $\nu$.  The non-negative integer $c_{\mu \nu}^{ \lambda}$ is
called the {\em \LR\ coefficient} of the partitions $\mu$, $\nu$, and
$\lambda$.  So $c_{\mu \nu}^{\lambda} = |LR(\mu, \nu;\lambda)|.$
\end{df}

\section {The Algorithm}
\begin{figure}[H]
\centering{\includegraphics[scale=0.3]{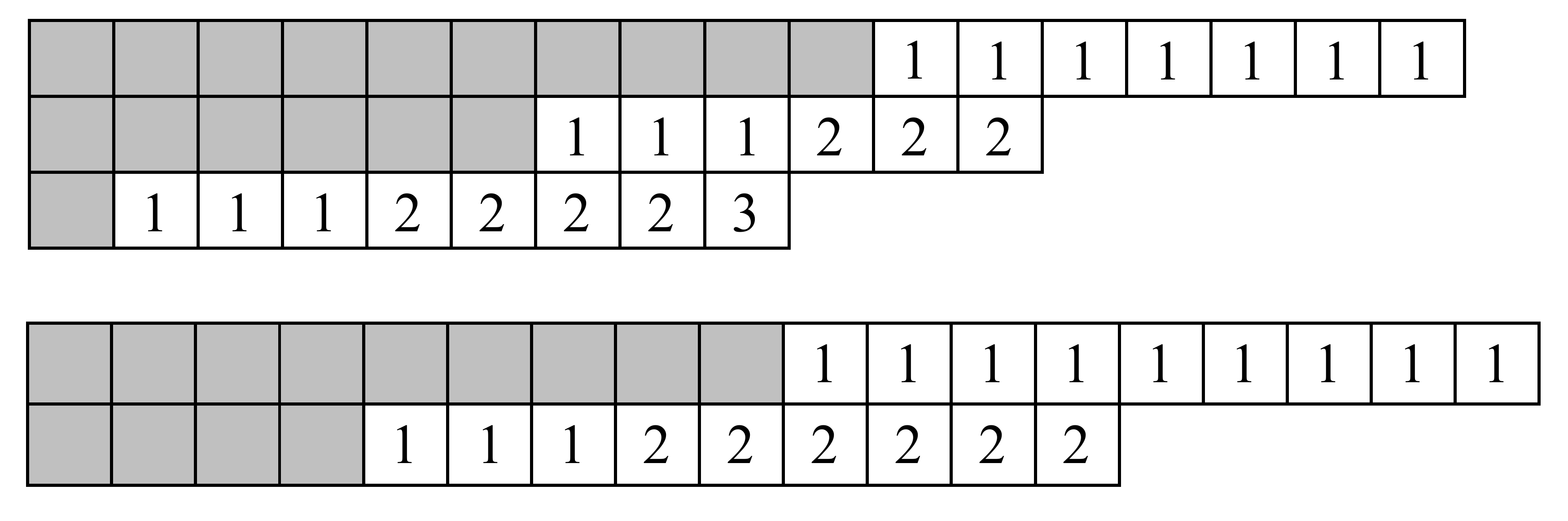}}
\caption{The two diagrams}
\label{Fig1}
\end{figure}

We now present the algorithm that takes two  \LR\ tableaux of types $(\mu, \nu;\lambda)$ and $(\mu', \nu';\lambda')$, and produces a \LR\ filling of type $(\mu\oplus \mu', \nu\oplus \nu' ;\lambda\oplus \lambda')$. We give the steps, along with an example for which
 $\mu=(10, 6, 1)$, $\nu=(13, 7, 1)$, $\lambda =(17, 12, 9)$
$\mu'=(9, 4)$, $\nu'=(12, 6)$, $\lambda' =(18, 13)$.  Hence
$\mu\oplus \mu'=(10, 9, 6, 4, 1)$, $\nu\oplus \nu' =(13, 12, 7, 6,1)$, and
$(\lambda\oplus \lambda')=(18, 17, 13, 12, 9)$.  The original diagrams are shown in Figure \ref{Fig1}.


\begin{figure}[H]
\centering{\includegraphics[scale=.3]{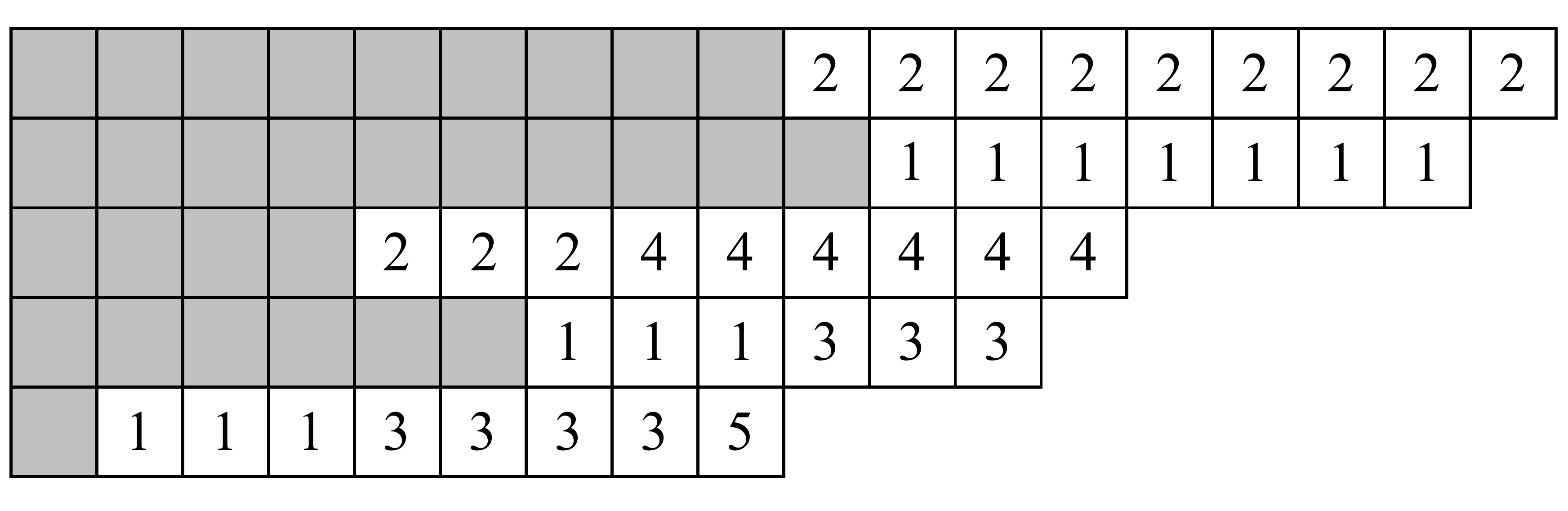}}
\caption{The initial summed diagram, with relabeled content}
\label{Fig2}
\end{figure}

The first step is to relabel the content of both diagrams so that the summed diagram will contain $1$'s, $2$'s, $\ldots$ and $n$'s where $n=$ length($\lambda$)+length($\lambda'$).  For our example, since $\nu_1 \geq \nu_{1}' \geq \nu_2 \geq \nu_{2}' \geq \nu_3$, the partition
$\nu$ gives us the number of $1$'s, $3$'s and $5$'s in the sum and $\nu'$ the number of $2$'s and $4$'s. The next step of the algorithm is to combine the parts of $\lambda$ and $\lambda'$ (without moving any entries within a given row) so that the row lengths weakly decrease. If two parts of $\nu$ and $\nu'$ are equal, we will choose the content so that the smaller number goes in a higher row of $\lambda \oplus \lambda'$. See Figure \ref{Fig2}.
Next, in order to make the parts formed by $\mu$ and $\mu'$ into the partition $\mu \oplus \mu'$, we sort the elements in each column of the diagram, putting the unfilled boxes of $\mu$ or $\mu'$ above any filled box, while retaining the order of the filled boxes within that column. Filled cells stay in their original columns, even if this leads to row-strict violations, as it does in rows 2 and 4 of Figure \ref{Fig3} below.
\begin{figure}[H]
\centering{\includegraphics[scale=.3]{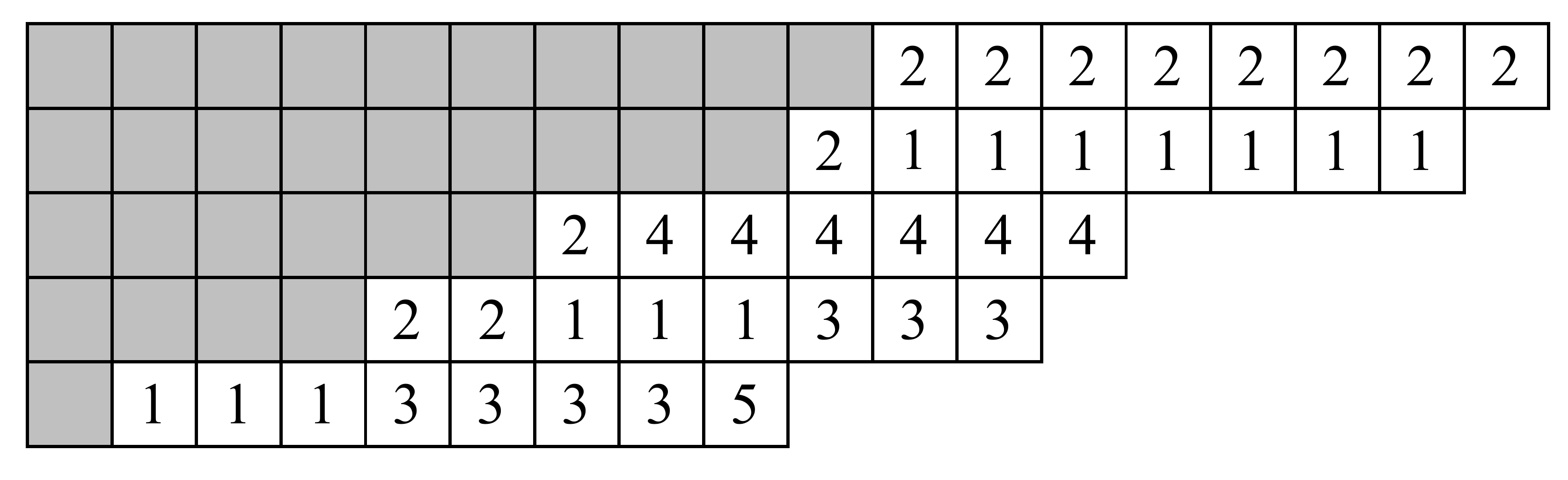}}
\caption{The summed diagram, after ordering the parts of $\mu$\label{Fig3}}
\end{figure}

At this point, the outer ($\lambda \oplus \lambda'$) and inner ($\mu \oplus \mu'$) shapes of the summed diagram are correct, but the filling can have word, column-strict and/or row-strict violations (corresponding to conditions (LR3), (LR2) and the containment requirement $\lambda^{(i)} \subseteq \lambda^{(i+1)}$, respectively, in the definition of \LR\ fillings).  We will first correct violations to the rule that entries in the boxes of row $i$ are less than or equal to $i$.  We correct these ``$i$ in row $i$" violations by following the ``north-east" rule:  Suppose there are some entries with value $(i+k)$ appearing in row $i$ (with $k>0$).  There cannot be more of these entries than there are $i$'s in the diagram in rows below row $i$.  We look for the north-east most box below row $i$ containing an $i$ such that the number of $i$'s weakly north-east of that box equals the number of $(i+k)$'s weakly north-east of it (which includes the $(i+k)$'s of ``$i$ in row $i$" violations in row $i$, but may also include $(i+k)$'s in lower rows). There is, therefore, a horizontal strip, or ``strand" of violating $(i+k)$ entries in row $i$ and below, and a horizontal strip of $i$'s of equal length. We then swap these strands.
In Figure \ref{Fig4} below, to fix the ``1 in row 1" violation, we will swap a strand of $1$'s and $2$'s.  We move the heavy frame from the last box of row 1 (starting in the first diagram) working south-west (into the second diagram), until in the third diagram the strands of $1$'s and $2$'s (the shaded boxes) have equal length.  Note that in our example we found the north-east most $1$ ending a strand of length equal to the strand of $2$'s above it. This strand of 2's includes a $2$ appearing in row 2.  Nevertheless, we then swapped the contents of the two strands as shown in Figure \ref{Fig5}, moving $2$'s out of row 1, as well as moving a $2$ out of row 2.
\begin{figure}[H]
\centering{\includegraphics[scale=.3]{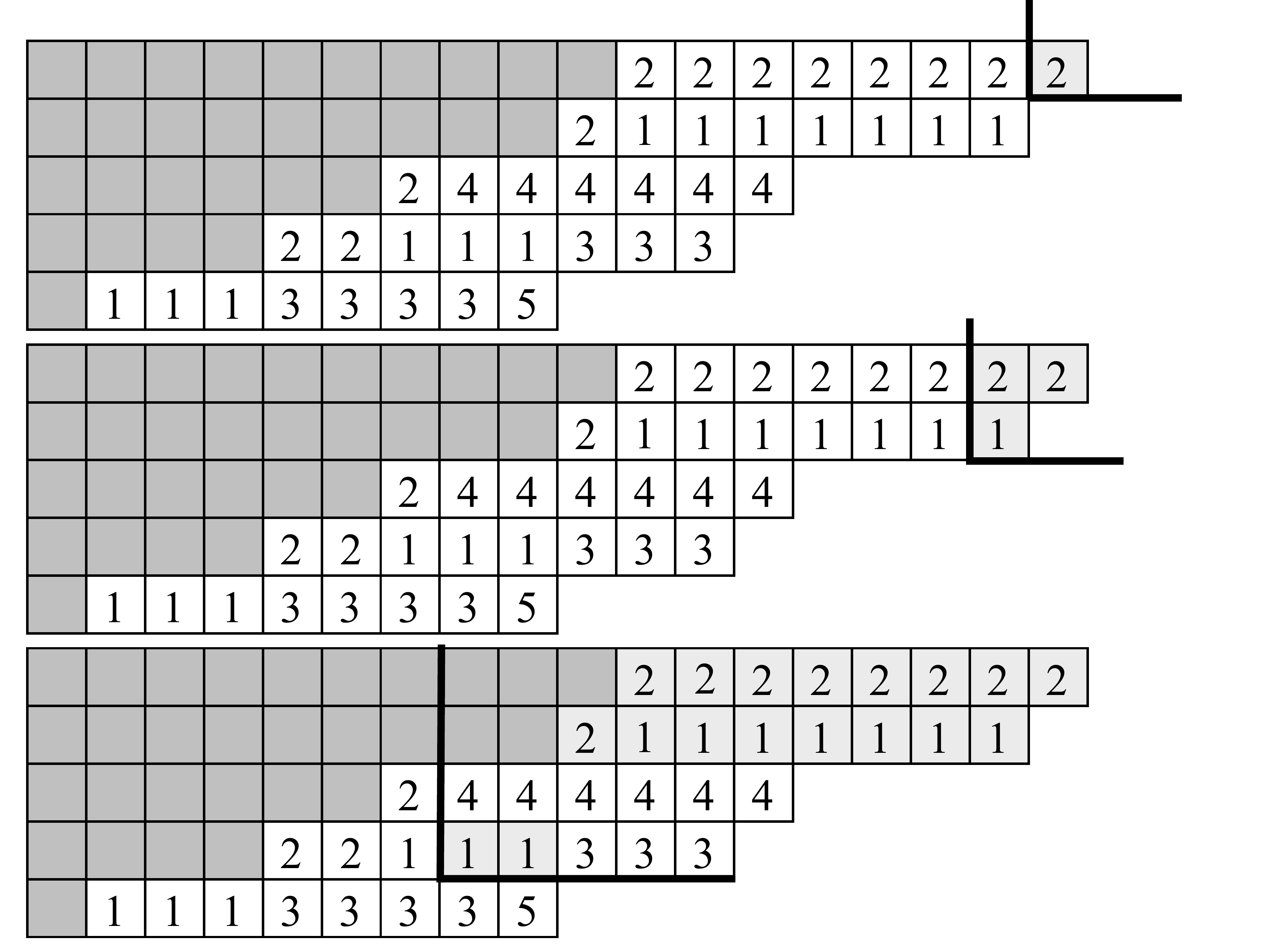}}
\caption{Sequence of summed diagrams, searching for equal strands of 1's and 2's (shaded)}
\label{Fig4}
\end{figure}
\begin{figure}[H]
\centering{\includegraphics[scale=.3]{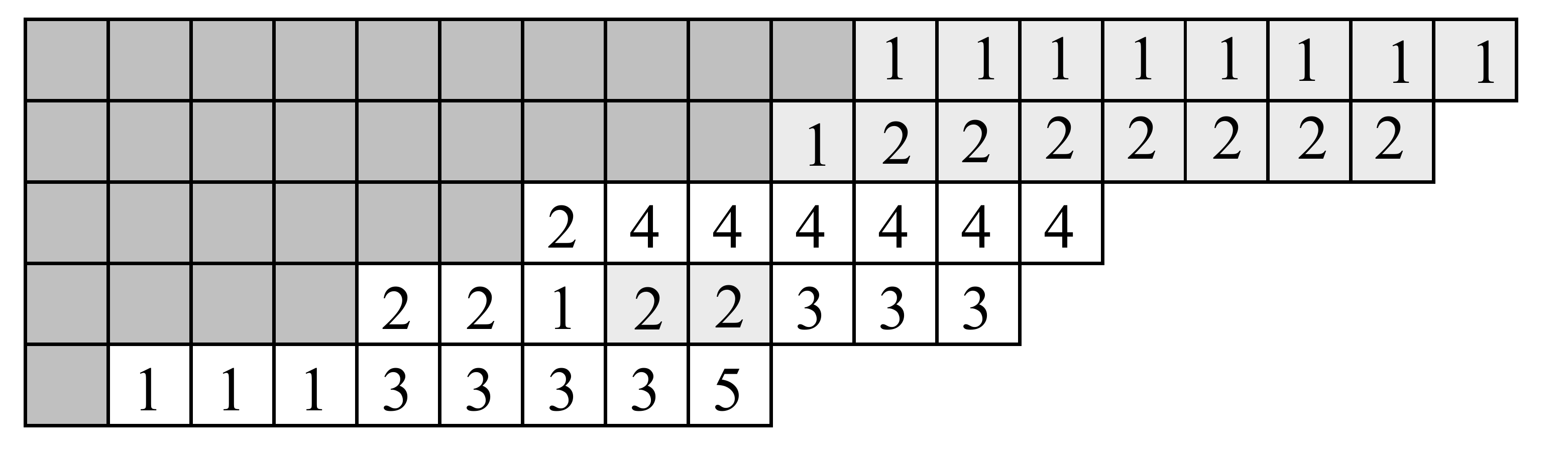}}
\caption{The summed diagram, with strands of 1's and 2's swapped}
\label{Fig5}
\end{figure}

Once we have replaced any $(i+k)$'s appearing in row $i$ with $i$'s, if $k-1>0$, we repeat this process, now replacing any $(i+k-1)$'s in row $i$ with a strand of $i$'s appearing in lower rows, etc., until there are only the numbers $1$ through $i$ appearing in row $i$. The feasibility of this (and every) step of the algorithm will be addressed in our proof below for the correctness of the algorithm.

\begin{figure}[H]
\centering{\includegraphics[scale=.3]{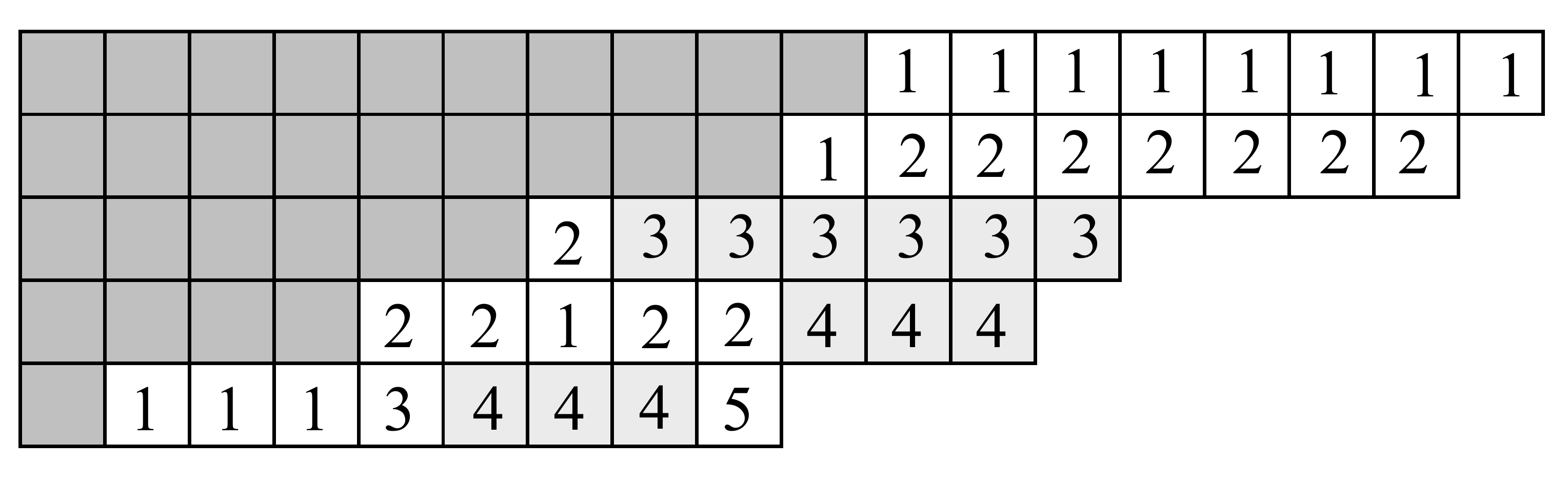}}
\caption{The summed diagram, with strands of 3's and 4's swapped}
\label{Fig6}
\end{figure}

We then continue this process, working down to the bottom row, until the content of each row is no bigger than the row index. In Figure \ref{Fig6}, we have swapped the $3$'s and $4$'s to fix the ``3 in row 3" problem.
The summed diagram may still contain column-strict, word, or row-strict violations. We will ignore the row-strict violations for now but will fix column-strict or word violations as they appear.  We shall do this by finding the north-east most ``bad box".  In the case of word violations, a bad box in row $j$, say, is one such that for some $i$, the number of $(i+1)$'s north-east of the bad box (including the bad box) exceeds the number of $i$'s through row $j-1$. For column strict violations a bad box is a box that contains an entry $i$ which lies directly above an entry $k$ where $k < i$.   If the north-east most bad box corresponds to a word violation, we will fix it by swapping strands as in the fixes for ``$i$ in row $i$" violations. In the case of column strict violations, we simply swap the two adjacent entries in the column so that they are strictly decreasing down the column. For example, in Figure \ref{Fig6}, working from the north-east corner, we would first encounter a
column-strict violation in the form of a pair of 3's above a pair of 2's, in rows 3 and 4, which we swap, shown in Figure \ref{Fig7}:
\begin{figure}[H]
\centering{\includegraphics[scale=.3]{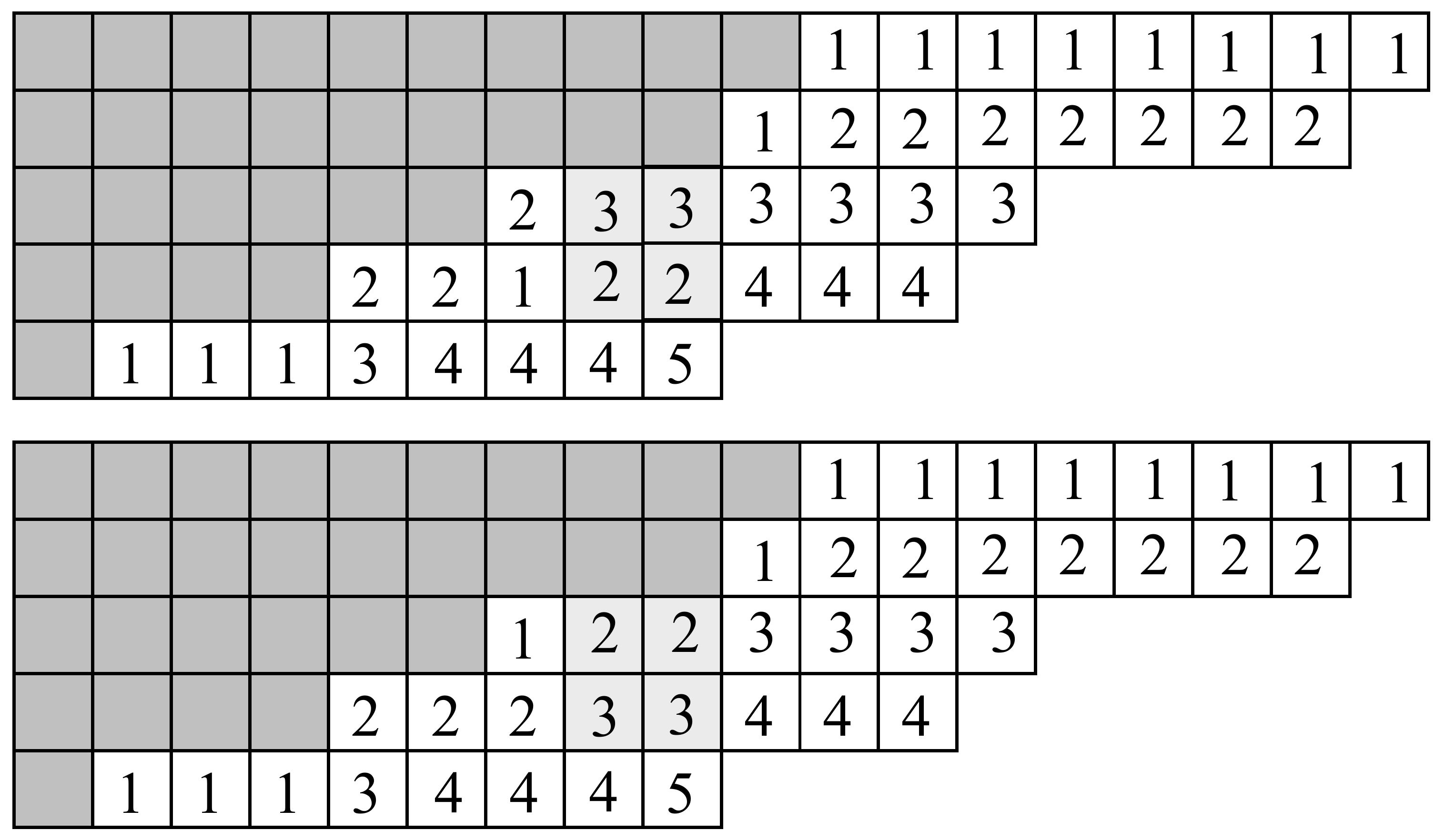}}
\caption{The summed diagram, with the pair of shaded 2's and 3's swapped}
\label{Fig7}
\end{figure}
We then encounter a word violation as the ten 2's in rows 2 and 3 exceed the nine 1's in rows 1 and 2.  This violation is caused by the first 2 in row 3, so begin moving the frame south-west from that ``bad" box, until the strands of 1's and 2's are equal in length.  The swap is shown in Figure \ref{Fig8} below:
\begin{figure}[H]
\centering{\includegraphics[scale=.3]{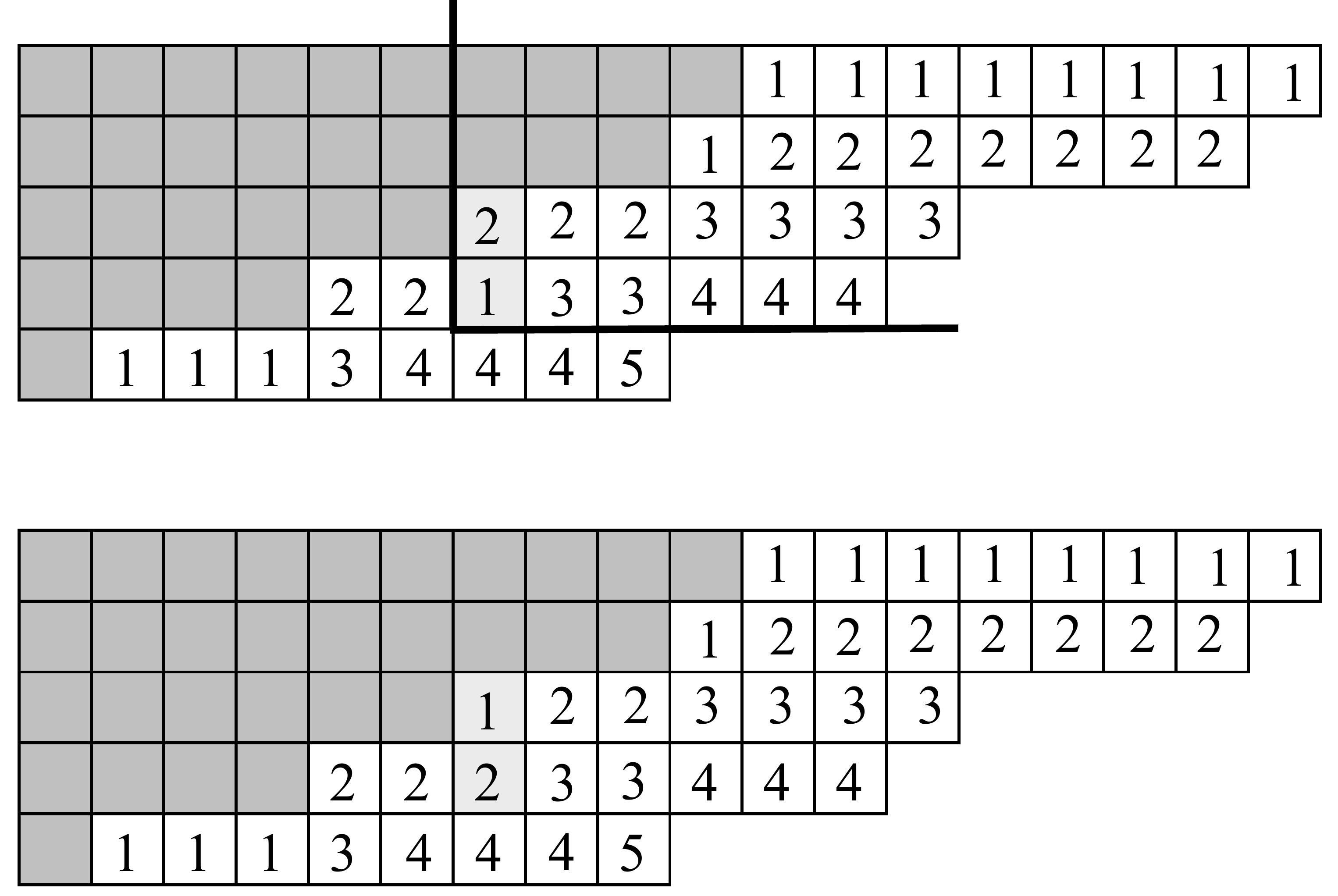}}
\caption{The summed diagram, with the shaded strands of 1's and 2's identified and swapped}
\label{Fig8}
\end{figure}
In our example, there is one final word violation (the 2's through row 4 exceed the 1's through row 3). The first bad box is the middle 2 in row 4, and we swap the shaded boxes, as shown in
Figure \ref{Fig9}.
\begin{figure}[H]
\centering{\includegraphics[scale=.3]{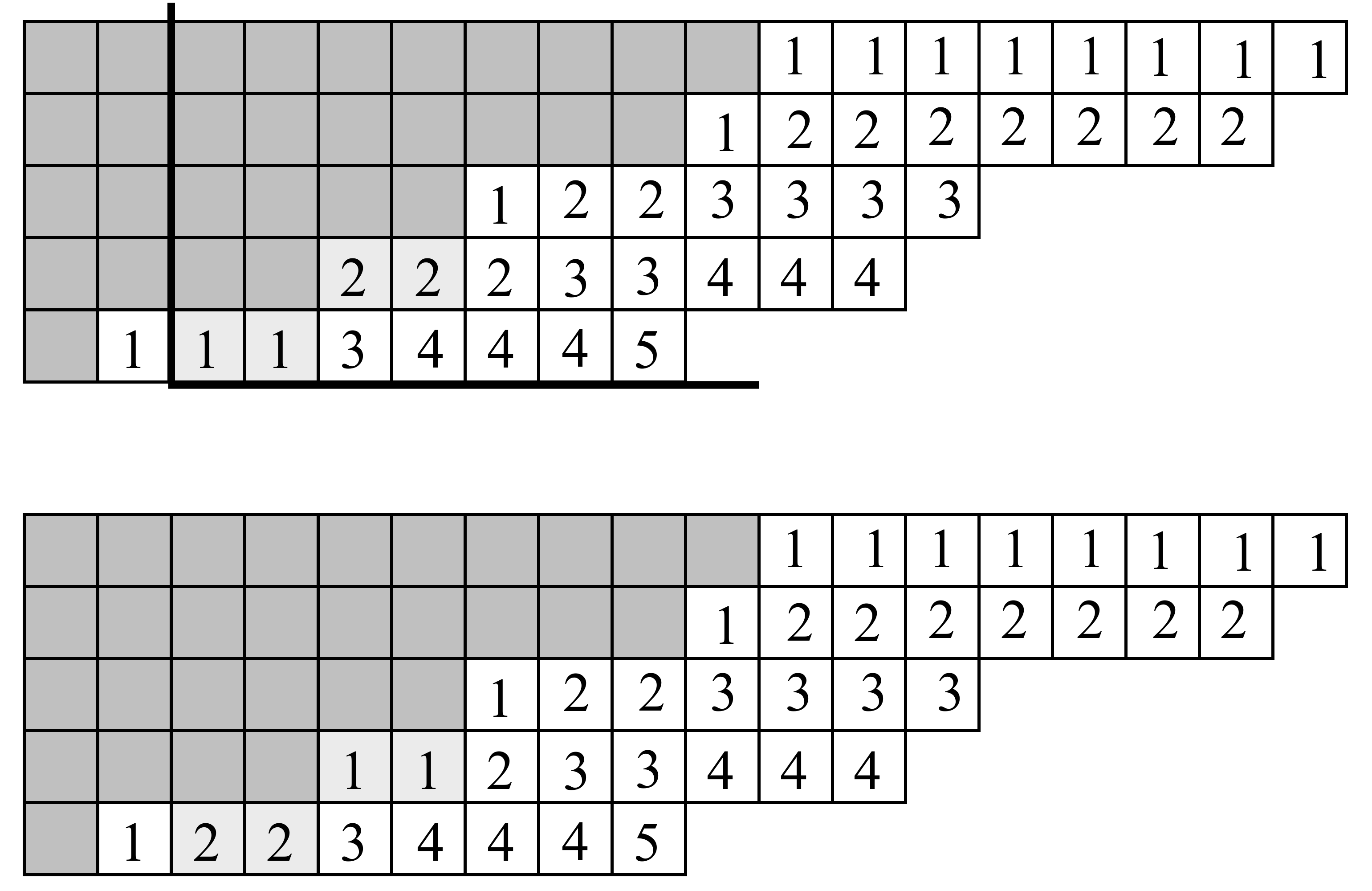}}
\caption{The summed diagram, with the final strands of 1's and 2's identified and swapped}
\label{Fig9}
\end{figure}
After fixing all word or column-strict violations, any remaining row-strict violations can be fixed, by ordering the contents of each row.  The algorithm terminates, and the summed filling, as we shall show, is \LR.

\section {The Dual Graph and Flows}
We will now proceed to develop the machinery necessary to prove that the above algorithm necessarily terminates in a \LR\ filling for the summed diagram.
In particular, we will review the construction of {\em hives}, first found in~\cite{knut} in the work of Knutson and Tao. (See also the work of Buch~\cite{Buch} for a more detailed discussion of hives, and Pak and Vallejo~\cite{pak} for their connection to \LR\ fillings.)

\begin{df} A {\em hive} of size $r$ is a triangular array of
numbers $(h_{ij})_{0 \leq i \leq j \leq r}$ that satisfy the {\em rhombus inequalities}:
\begin{enumerate}
\item {\em Right:} $h_{ij}-h_{i,j-1} \geq h_{i-1, j} - h_{i-1,j-1}$, for $1
\leq i < j \leq r$.
\item {\em Vertical:} $h_{i-1,j} - h_{i-1,j-1} \geq h_{i,j+1} - h_{ih}$, for $1
\leq i \leq j \leq r$.
\item {\em Left:} $h_{ij} - h_{i-1,j} \geq h_{i+1,j+1} - h_{i,j+1}$, for $1
\leq i \leq j < r$.
\end{enumerate}
\end{df}
A hive of size 5 is shown below.
\[ \setlength{\arraycolsep}{-0.05mm} \renewcommand{\arraystretch}{1.3}
\begin{array}{ccccccccccc}
    &     &    &    &     &h_{00}&     &     &      &   &  \\
    &     &    &    &h_{10}&     &h_{11}&     &      &   &  \\
    &     &    &h_{20}&     &h_{21}&     &h_{22}&      &   &  \\
    &     &h_{30}&    &h_{31}&     &h_{32}&     &h_{33}&   &  \\
    &h_{40}&    &h_{41}&     &h_{42}&     &h_{43}&      &h_{44}&  \\
h_{50} & & h_{51} & & h_{52} & & h_{53} & & h_{54} && h_{55} \end{array} \]

Note that we will always normalize hives so that $h_{00}=0$.

These rhombus inequalities are so named because the terms in each inequality form a rhombus in the hive, made by adjacent entries in the array such that that the upper acute angle points to the right,
vertically, and to the left, respectively.  In each case the inequality asserts
that the sum of the entries of the obtuse vertices of the rhombus is greater
than or equal to the sum of the acute entries.

Note that the definition of a hive allows its entries to be real numbers, but we shall restrict our attention to hives whose entries are nonnegative integers.

Let $HIVE_{r}$ denote the set of nonnegative integer-valued hives of size
$r$.  We shall say a hive $H$ is of {\em type} $\tau(H)= (\mu,\nu,\lambda)$ (for
sequences of nonnegative integers $\mu$, $\nu$, and $\lambda$ of length $r$)
when
\begin{align*}
\mu_i &= h_{i0} - h_{(i-1),0} \quad & \hbox{(the downward
differences of entries
along the left side)} \\
\nu_i &= h_{r i} - h_{r(i-1)}, \quad & \hbox{(the rightward
differences of entries along the bottom)}\\ \lambda_{i} &=
h_{(i)(i)} - h_{(i-1)( i-1)} \quad & \hbox{(the downward differences of
entries along the right side)}.
\end{align*}

Let
\[ H(\mu, \nu, \lambda)_r = \{ H \in HIVE_{r} : \tau(H) = (\mu, \nu,\lambda) \}. \]

As a consequence of the rhombus inequalities,  the type of a hive $H$ will be a triple of partitions of {\em weakly decreasing} nonnegative integers and, thus, form a triple of partitions.  A
necessary condition for the existence of a hive $H \in HIVE(\mu,\nu,\lambda)_r$ is:
\[ |\mu | + |\nu| = | \lambda|. \]

As an example, the hive below is an element of $HIVE(\mu,\nu,\lambda)_5$
where $\mu=(10,9,5,3,1)$, $\nu=(12, 11, 7, 6,1)$ and $\lambda=(18, 16, 12, 11, 8)$.
\vspace{-0.275in}\begin{figure}[H]
\[ \setlength{\arraycolsep}{-0.05mm} \renewcommand{\arraystretch}{1.3}
H = \begin{array}[h]{ccccccccccc}
    &     &    &    &     &0&     &     &      &   &  \\
    &     &    &    &10&     &18&     &      &   &  \\
    &     &    &19&     &27&     &34&      &   &  \\
    &     &24&    &34&     &42&     &46&   &  \\
    &27&    &38&     &48&     &54&      &57&  \\
28 & & 40 & & 51 & & 58 & &64 && 65
 \end{array} .\] \caption{An example of a hive in $HIVE_{5}$}\label{hive-ex}
 \end{figure}

Pak and Vallejo~\cite{pak} gave an injective map $\phi$ from $LR(\mu, \nu;\lambda)$  to $H(\mu, \nu, \lambda)_r$ with
$$\phi(\left\{k_{ij}\right\}) =\left\{h_{ij}\right\},$$
defined by
\begin{equation}\label{h_pq}
h_{pq}=\sum_{i=1}^q\sum_{j=1}^p {k_{ij}}+\sum_{s=1}^p \mu_s.
\end{equation}
For example,

\bigskip
\hspace{1.8in} \vbox{ \offinterlineskip \openup-0.5pt
\nl
\f\f\f\f\f\f\f\f\f\f\b{$1$}\b{$1$}\b{$1$}\b{$1$}\b{$1$}\b{$1$}\b{$1$}\b{$1$}\nl
\f\f\f\f\f\f\f\f\f\b{$2$}\b{$2$}\b{$2$}\b{$2$}\b{$2$}\b{$2$}\b{$2$}\nl
\f\f\f\f\b{$1$}\b{$1$}\b{$1$}\b{$2$}\b{$3$}\b{$3$}\b{$3$}\b{$3$}\nl
\f\f\f\b{$1$}\b{$2$}\b{$2$}\b{$3$}\b{$3$}\b{$4$}\b{$4$}\b{$4$}\nl
\f\b{$1$}\b{$2$}\b{$3$}\b{$4$}\b{$4$}\b{$4$}\b{$5$}}

\bigskip\noindent
is the \LR\ tableau of type $((10, 9, 5, 3, 1), (12, 11, 7, 6, 1); (18, 16, 12, 11, 8))$ corresponding
to the hive given above.  Note that $h_{pq}$ equals the sum of the parts of $\mu$ through row $p$, and the $1$'s,
$2$'s, $\dots$, and $q$'s in rows 1 through $p$.  So for example, $h_{52}$ equals
$\mu_1+\dots \mu_5 +$ (the total number of $1$'s and $2$'s through row $5$).

In fact, $\phi$ is onto the set of non-negative integer-valued hives, and so we also have
$$\phi^{-1}(\{h_{pq}\})=\{k_{ij}\}$$
where
\begin{equation}\label{k_ij}
k_{ij}=(h_{(j-1)(i-1)}+h_{j i})-(h_{(j-1) i}+h_{j(i-1)})
\end{equation}
for $i<j$.  This difference is in fact the rhombus difference for right-slanted rhombi, and so is non-negative.  In the example above, $k_{24}=(h_{31}+h_{42})-(h_{32}+h_{41})=(34+48)-(42+38)=2$

In order to verify our algorithm does indeed terminate in a \LR\ filling, we will need several more definitions, which will re-formulate hive combinatorics in a more convenient form.
We can view a hive as a vertex labeling on an underlying undirected graph, as shown:
\begin{figure}[H]
\centering{\includegraphics[scale=.3]{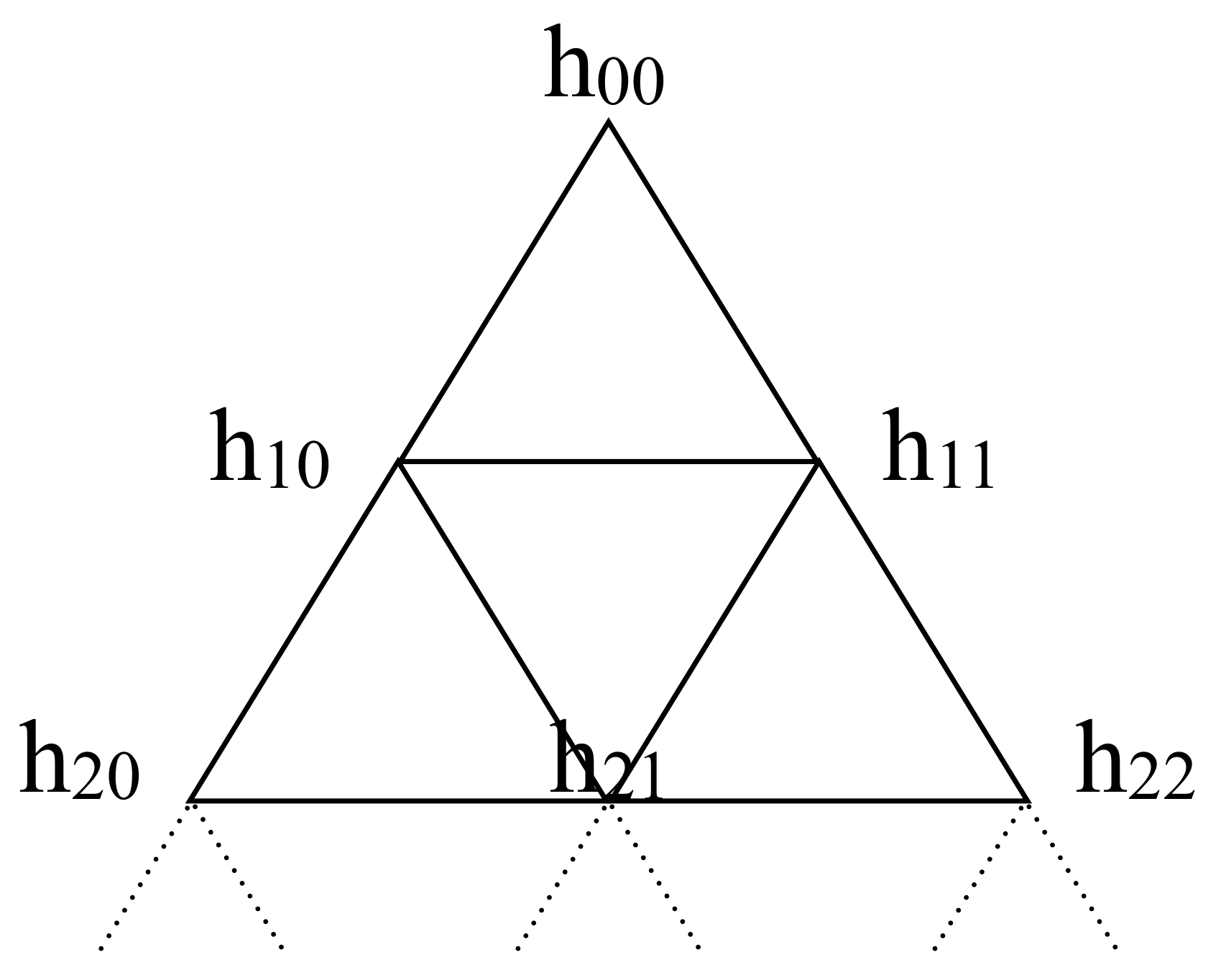}}
\caption{The hive graph}
\label{Fig10}
\end{figure}
So, for example, our hive shown in Figure~\ref{hive-ex} corresponds to the hive graph:
\begin{figure}[H]\
\centering{\includegraphics[scale=.2]{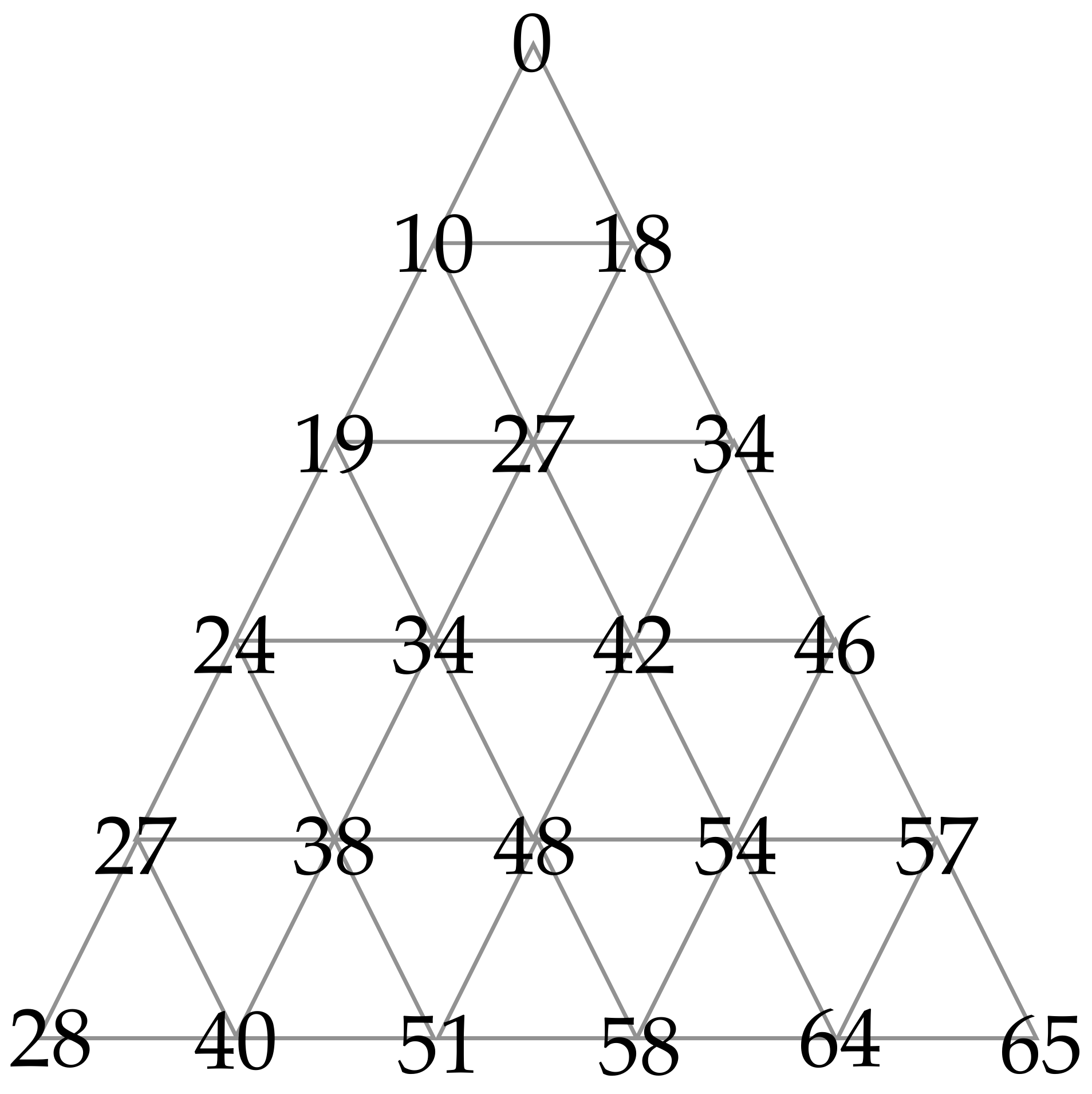}}
\caption{The hive graph for the hive in Figure~\ref{hive-ex}}
\label{Fig10a}
\end{figure}
We then can create the dual graph to the hive graph:

\begin{figure}[H]
\centering{\includegraphics[scale=.25]{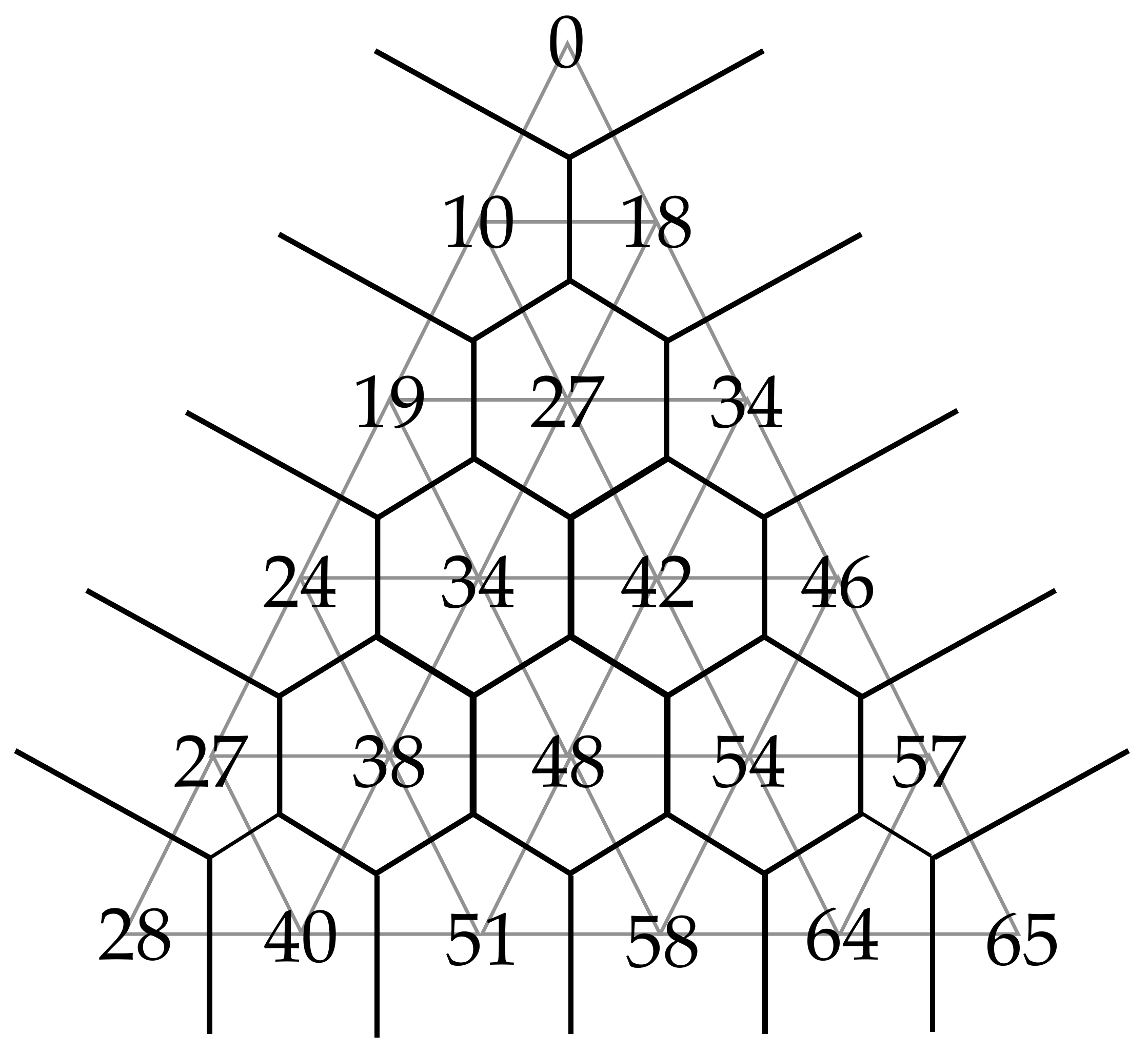}}
\caption{The dual graph}
\label{Fig11}
\end{figure}

Each hive value now lies in a unique cell of the dual graph. We weight each edge in the dual graph by the positive difference of the adjacent hive entries from the cells on either side of the edge:
\begin{figure}[H]
\centering{\includegraphics[scale=.3]{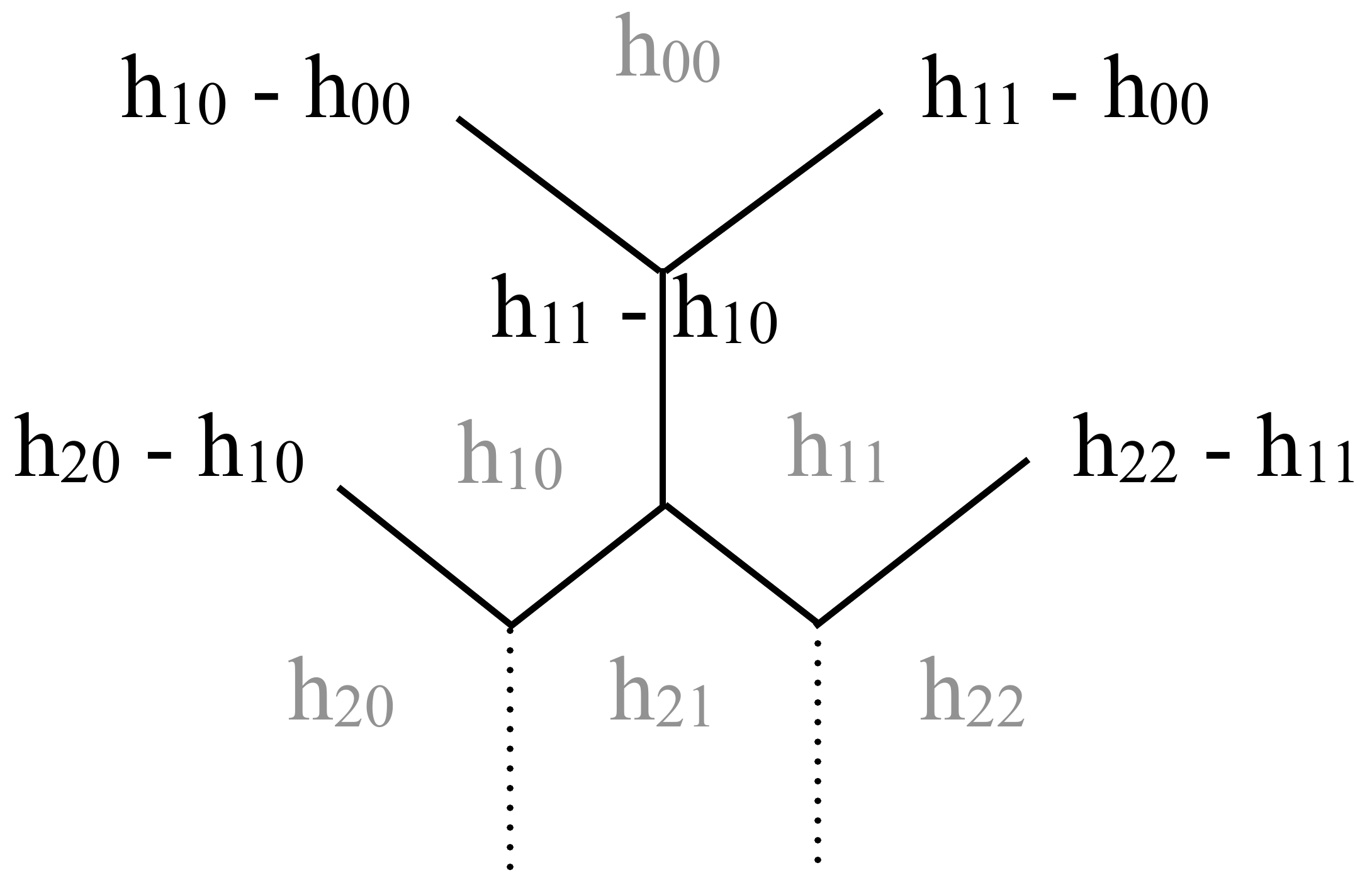}}
\caption{The weighted dual graph}
\label{Fig12a}
\end{figure}
In the diagram above, the grey numbers are the hive entries, the dark numbers are the corresponding edge weights of the dual graph.

So for our hive example, the top portion of the weighted dual graph is:
\begin{figure}[H]
\centering{\includegraphics[scale=.2]{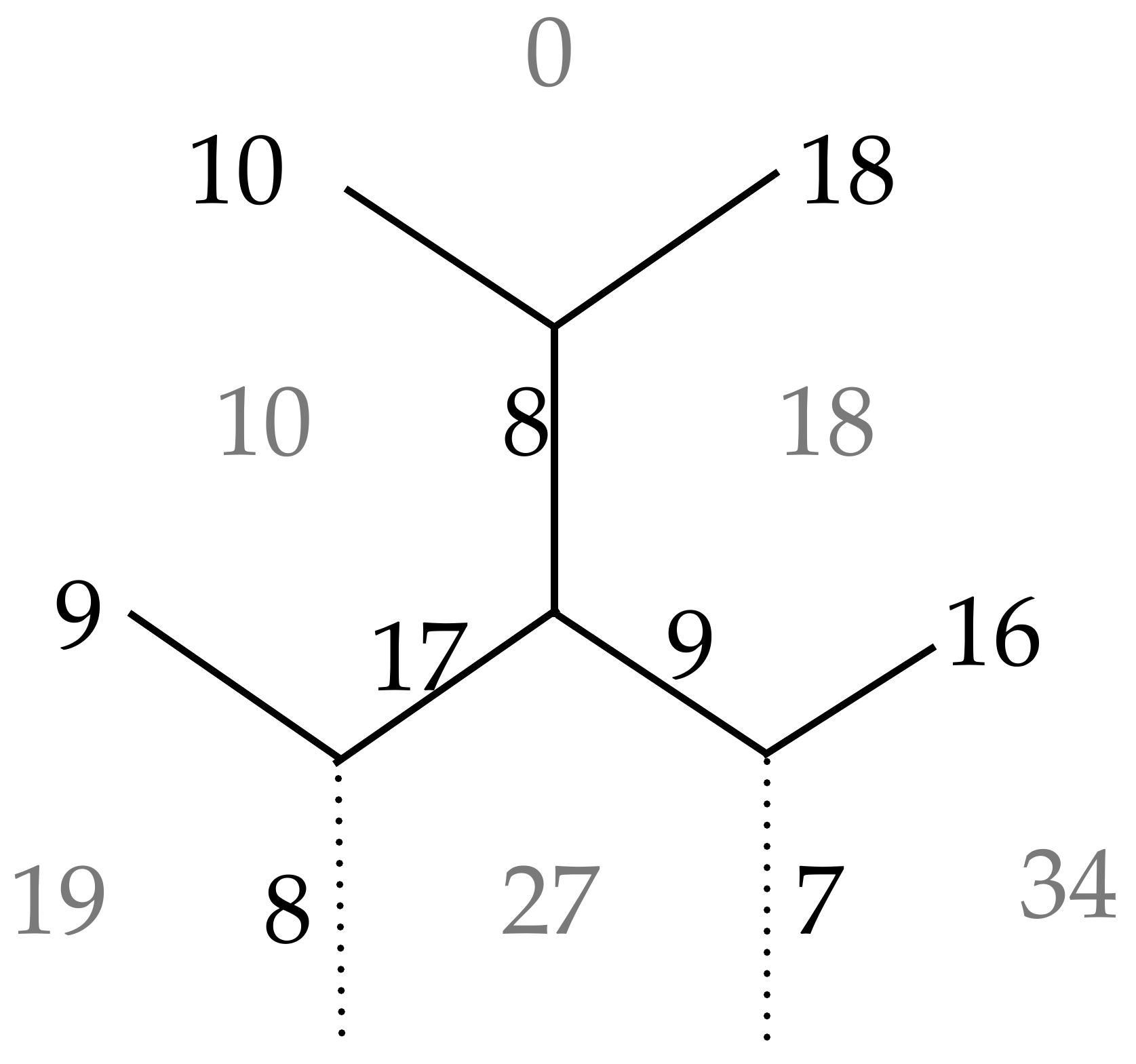}}
\caption{The top of the weighted dual graph}
\label{Fig12}
\end{figure}
By using Equation \ref{h_pq}, we can compute, in terms of the filling, the weights of the edges on the exterior of the dual graph. The weights of the left edges are the parts of $\mu$, the bottom weights are the parts of $\nu$ and the right edge weights are the parts of $\lambda$.  Similarly, we may find formulas for edge weights for the interior edges. Using the diagram below to identify certain edge weights in a typical dual graph cell,
\begin{figure}[H]
\centering{\includegraphics[scale=.25]{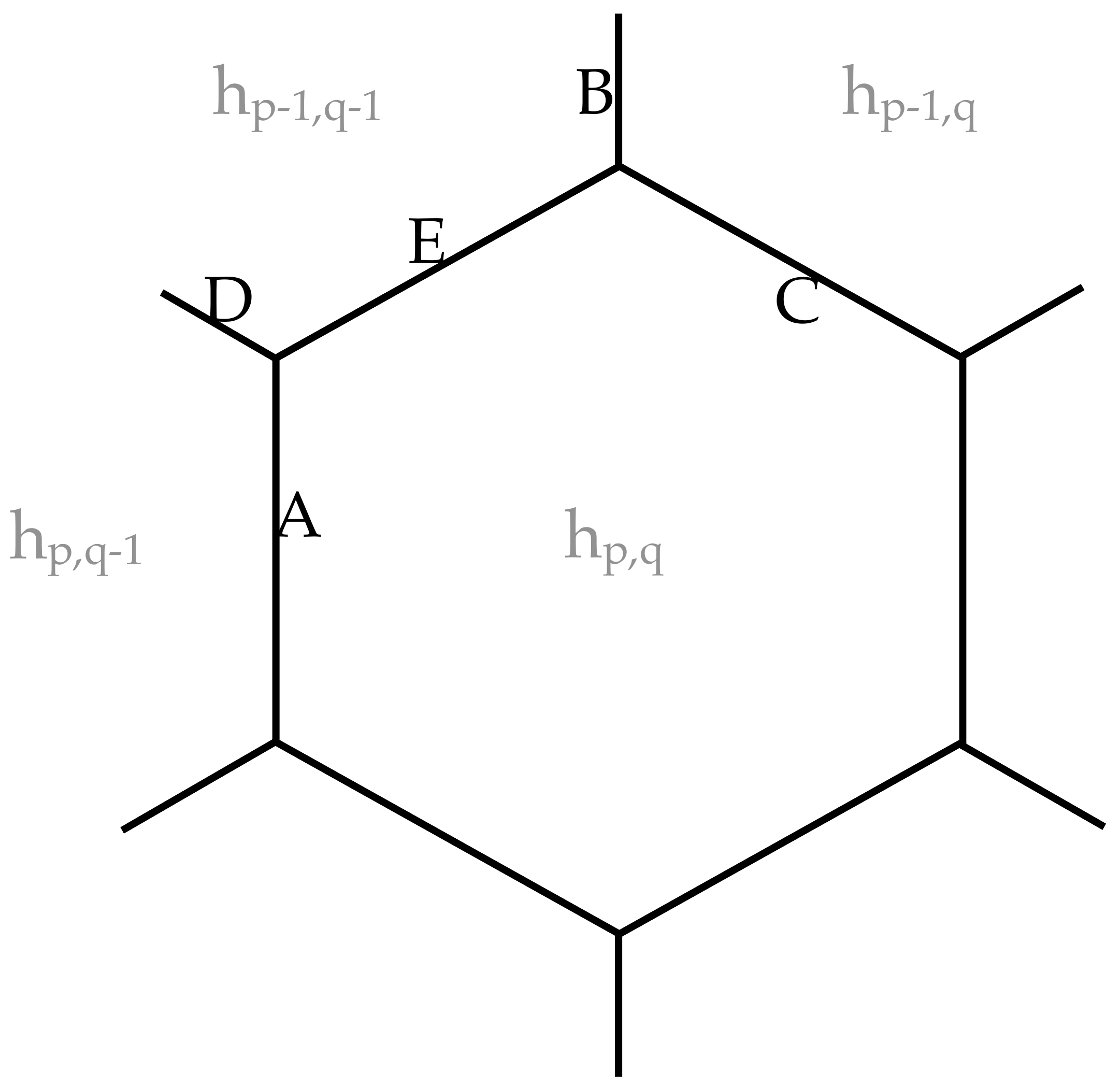}}
\caption{Edge weights identified}
\label{Fig23}
\end{figure}
\noindent we have (from Equation \ref{k_ij}), that, for $q < p$, $k_{qp} = A - B = C - D$.  Then, we have formulas for the edge weights $A$, $B$, $C$ and $E$:
\begin{eqnarray}
A &= &k_{q,q}+k_{q, q+1}+\cdots+k_{q,p}\\
B &= &k_{q,q}+k_{q, q+1}+\cdots+k_{q,p-1}\\
C &=& \mu_p + k_{1,p}+k_{2,p}+\cdots+k_{q,p}\\
E &=& \mu_p + k_{1,p}+k_{2,p}+\cdots+k_{q-1,p}+k_{q,q}+k_{q,q+1}+\cdots k_{q,p}
\end{eqnarray}

We can view each edge weight as a capacity for a flow, with edges flowing in the directions as shown:
\begin{figure}[H]
\centering{\includegraphics[scale=.3]{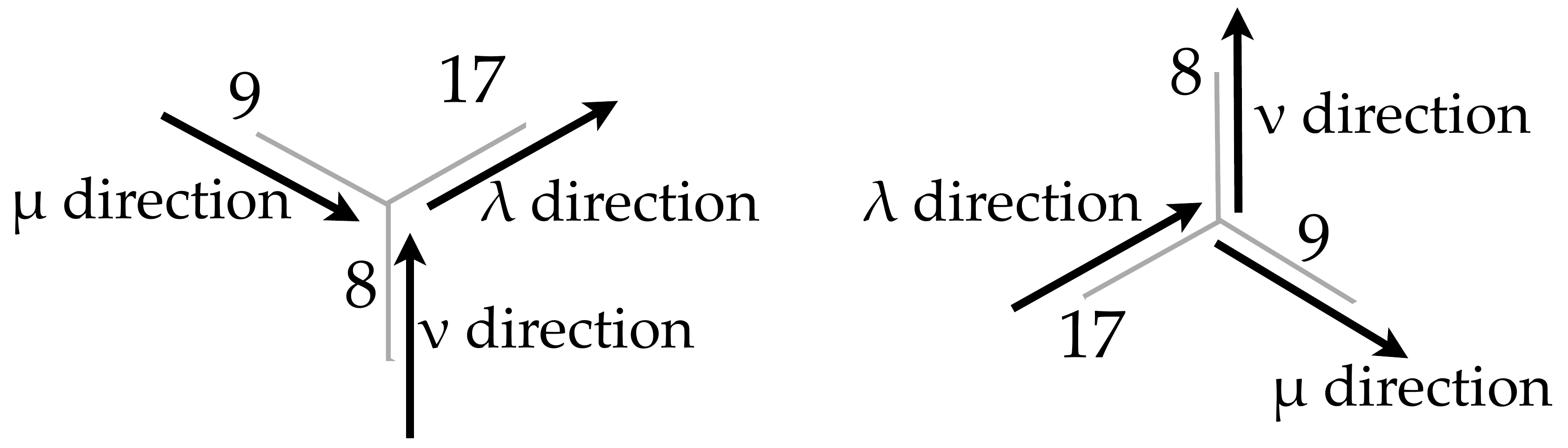}}
\caption{Dual graph edges, with capacities and flow directions labeled}
\label{Fig13}
\end{figure}
Note that, because the edge capacities are defined as differences of the adjacent hive entries, the sum of the capacities flowing {\em into} a node equals the capacity {\em out} of the node.

For our hive example, the dual graph with flow capacities is given below:
\begin{figure}[H]
\centering{\includegraphics[scale=.3]{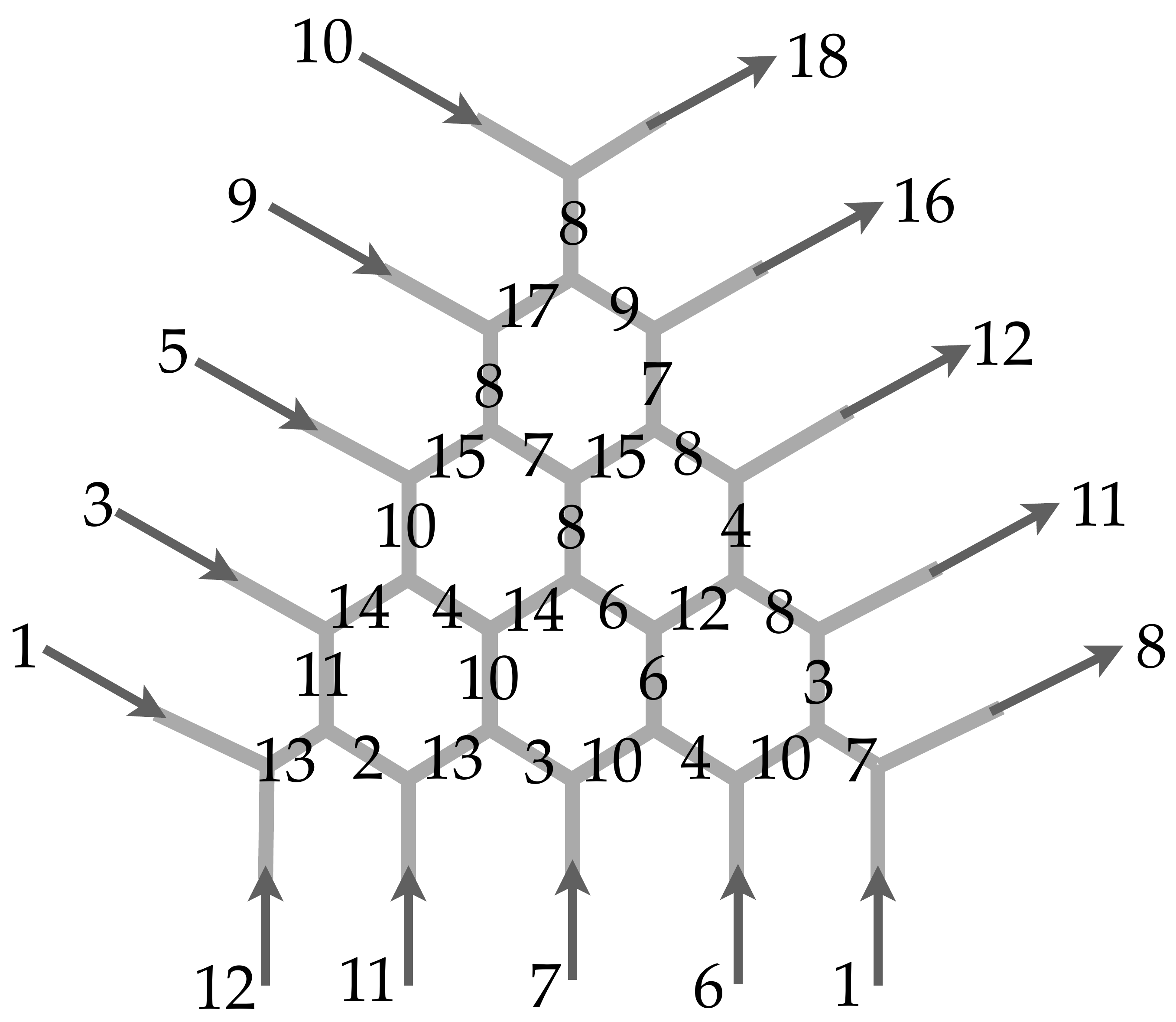}}
\caption{Dual graph edges, with flows labeled}
\label{Fig14}
\end{figure}

\begin{figure}[H]
\centering{\includegraphics[scale=.4]{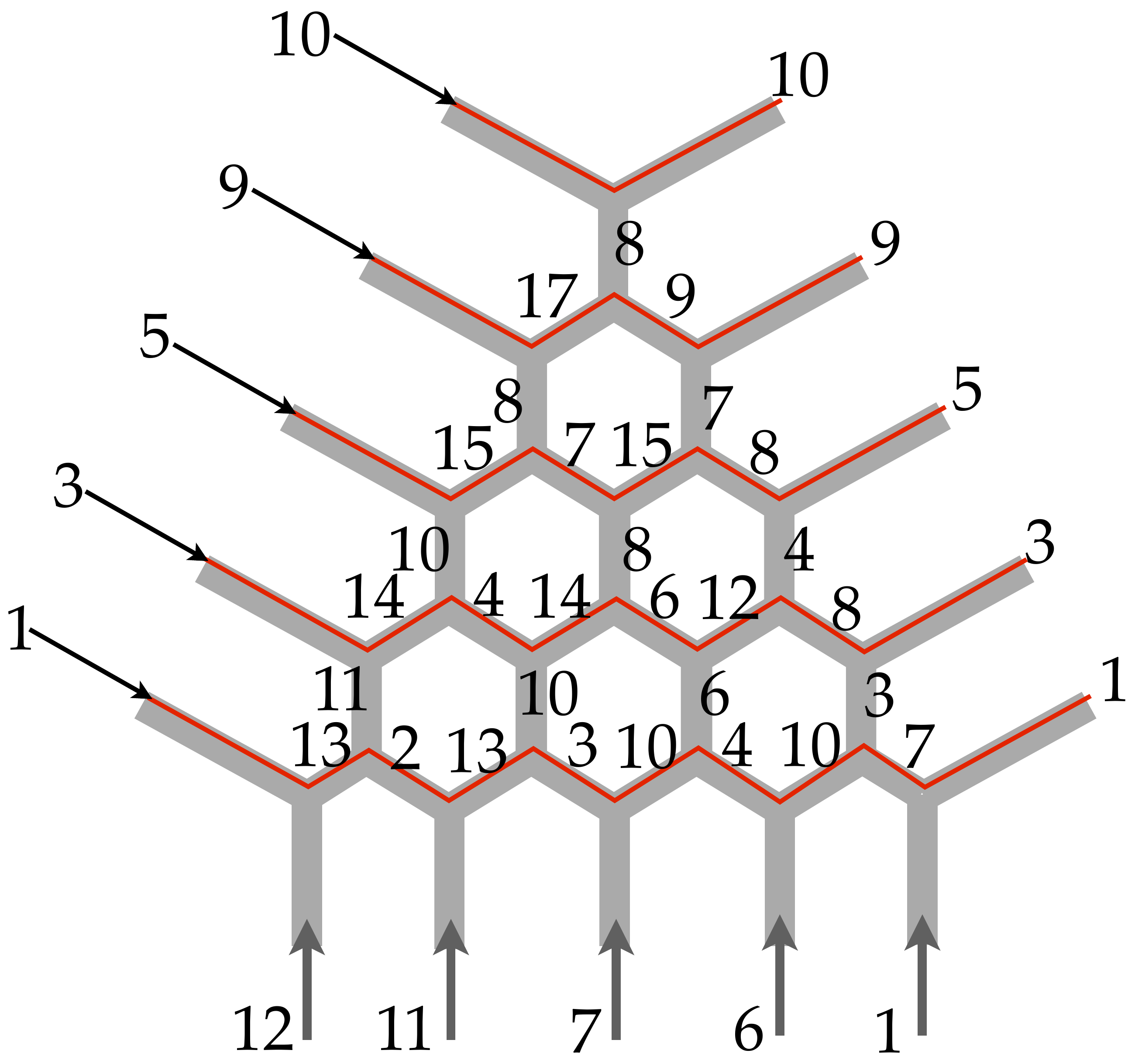}}
\caption{$\mu$ flow through the dual graph}
\label{Fig15}
\end{figure}

We can use Equations 3 through 6 to assign to each edge of the weighted dual graph a decomposition of the flow on that edge determined by the underlying \LR\ filling. First, the flow coming in from the $\mu$ direction will be required to flow across the ``rows"  of the dual graph, so that the flow for $\mu_i$ flows out of $\lambda_i$.  Essentially, the flow changes direction (southeast or northeast) at each vertex. This contributes to the flow out of each part of $\lambda$.  We will call this the ``$\mu$ flow".  Equations 5 and 6 above demonstrate that the capacity across the row is large enough for the $\mu$ flow. See Figure \ref{Fig15}.
\medskip

The flow from the parts of $\nu$ canonically flows up the ``spines" of the dual graph, which is a path starting a some bottom vertical edge with initial weight $\nu_i$, and then proceeding by going north, and then north-east alternately at each vertex. However, we can encounter vertices where, if the flow from the $\nu$ direction continues up, we would exceed the capacity of the next edge.  In that case, the flow bifurcates as shown in the example below.

\begin{figure}[H]
\centering{\includegraphics[scale=.3]{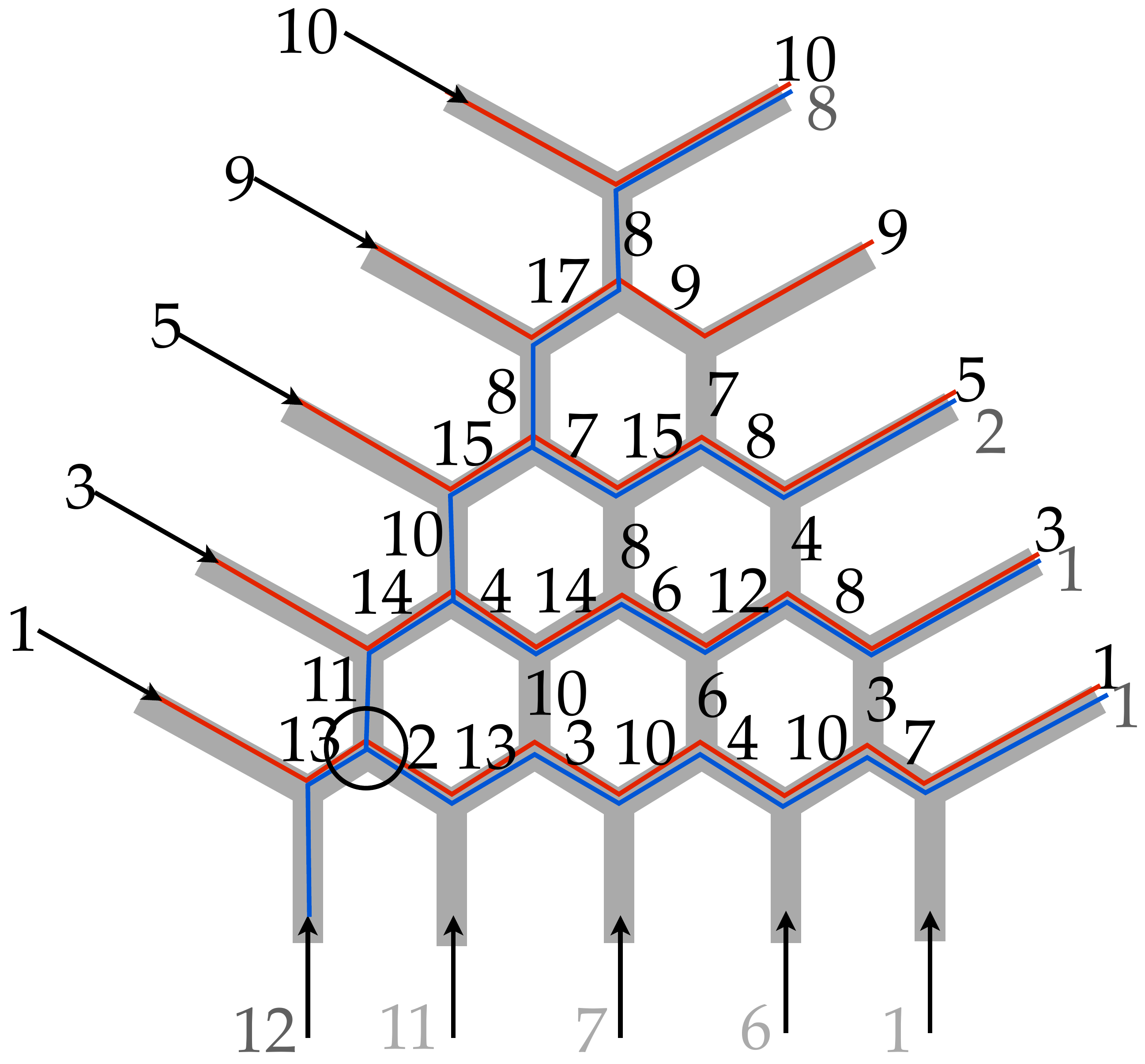}}
\caption{The $1$ flow through the dual graph}
\label{Fig16}
\end{figure}

Consider the flow coming in from $\nu_1$. In our example, $\nu_1 = 12$, and this corresponds to the twelve $1$'s in the associated \LR\ diagram.  The $1$'s should flow up the spine of the dual graph, but as shown in Figure \ref{Fig16} at the circled vertex, the flow of twelve $1$'s is more than the capacity of the vertical edge labeled with a capacity of 11.  Hence, the flow of $1$'s bifurcates, with the excess flow of one 1 going right, following the $\mu$ flow, and the remaining eleven $1$s continuing up the spine.  Note that the amount of flow that exceeds capacity of the following edge can always bifurcate in this manner because the capacity into any vertex equals the capacity out.  And, again, the edge weights given determined by Equations 2 through 6 demonstrate that this bifurcation is mirrored in the \LR\ filling.  For the 1's, at the circled vertex, the vertical edge weight of the dual graph has capacity $k_{11}+k_{1,2}+\cdots +k_{1,p}$ (where $p$ is the number of rows in $\lambda\oplus\lambda')$ and the edge to the right has weight
$\mu_p+k_{1,p}$ (here $q=1$), which equals the flow on that edge.  The bifurcation of the $1$ flow is shown below, in blue. At each vertex along the left-most spine, the flow of $1$ bifurcates as needed, with the number of $1$'s that exceed capacity of the next vertical edge being diverted to the right, and following the $\mu$ flow of that row. Comparing Equations 3 and 4 shows that the amount of the bifurcation coming from $\nu_{i}$ in the $j$th row is precisely $k_{ij}$, and Equations 5 and 6 show that there is sufficient capacity in the horizontal flows at each bifurcation to carry this part of the flow.
In Figure~\ref{Fig16}, the number of $1$'s flowing out of a part  $\lambda_i$ is indicated by the second number
at that outflow.
Hence, eight $1$'s flow out of $\lambda_1$, no $1$'s flow out of $\lambda_2$, two $1$'s flow out of $\lambda_3$, and both $\lambda_4$ and $\lambda_5$ have one 1.  This agrees with the \LR\
filling.

We create the flow for the remaining parts of $\nu$ (the $2$'s through $5$'s, or in general, the $2$'s through the number of parts of $\nu$) in the same manner.  For example, the flow of the $2$'s (the green lines) is included in the diagram below.
\begin{figure}[H]
\centering{\includegraphics[scale=.3]{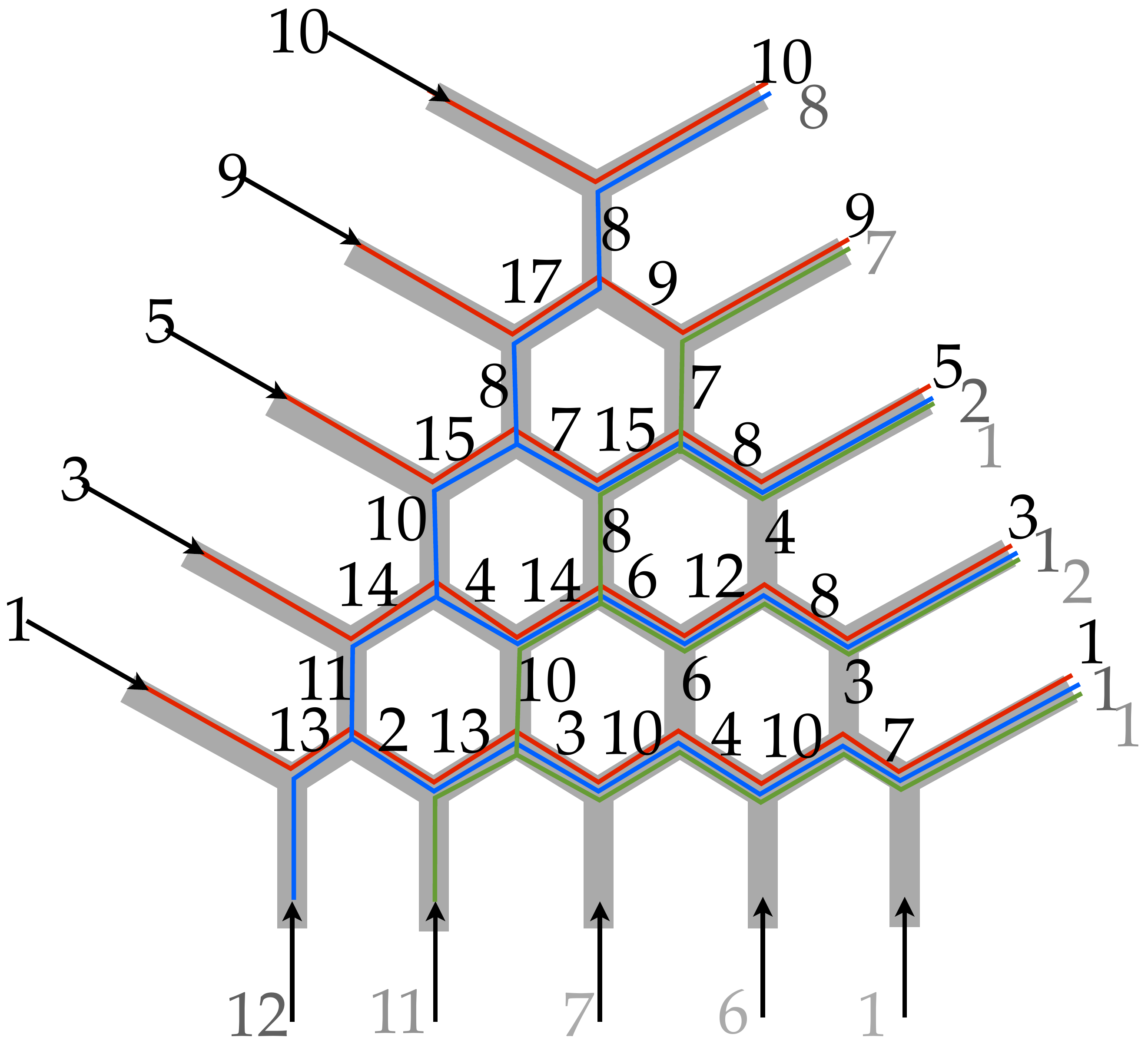}}
\caption{The $1$ and $2$ flow through the dual graph}
\label{Fig17}
\end{figure}
When the flow for all parts of $\nu$ has been allocated, we call this the ``canonical flow" of the dual graph;  it is uniquely determined by the \LR\ filling by means of Equations 2 through 6.   It is easy to show that the canonical flow of a weighted dual graph (subject to the inflow/outflow constraints at vertices) produces a \LR\ filling as well.  Burgisser and Ikenmeyer ~\cite{BH} used flows on the dual graph to analyze the positivity of \LR\ coefficients.

\section {From Flows on Dual Graphs to Flows on Honeycombs}
\begin{figure}[H]
\centering{\includegraphics[scale=.2]{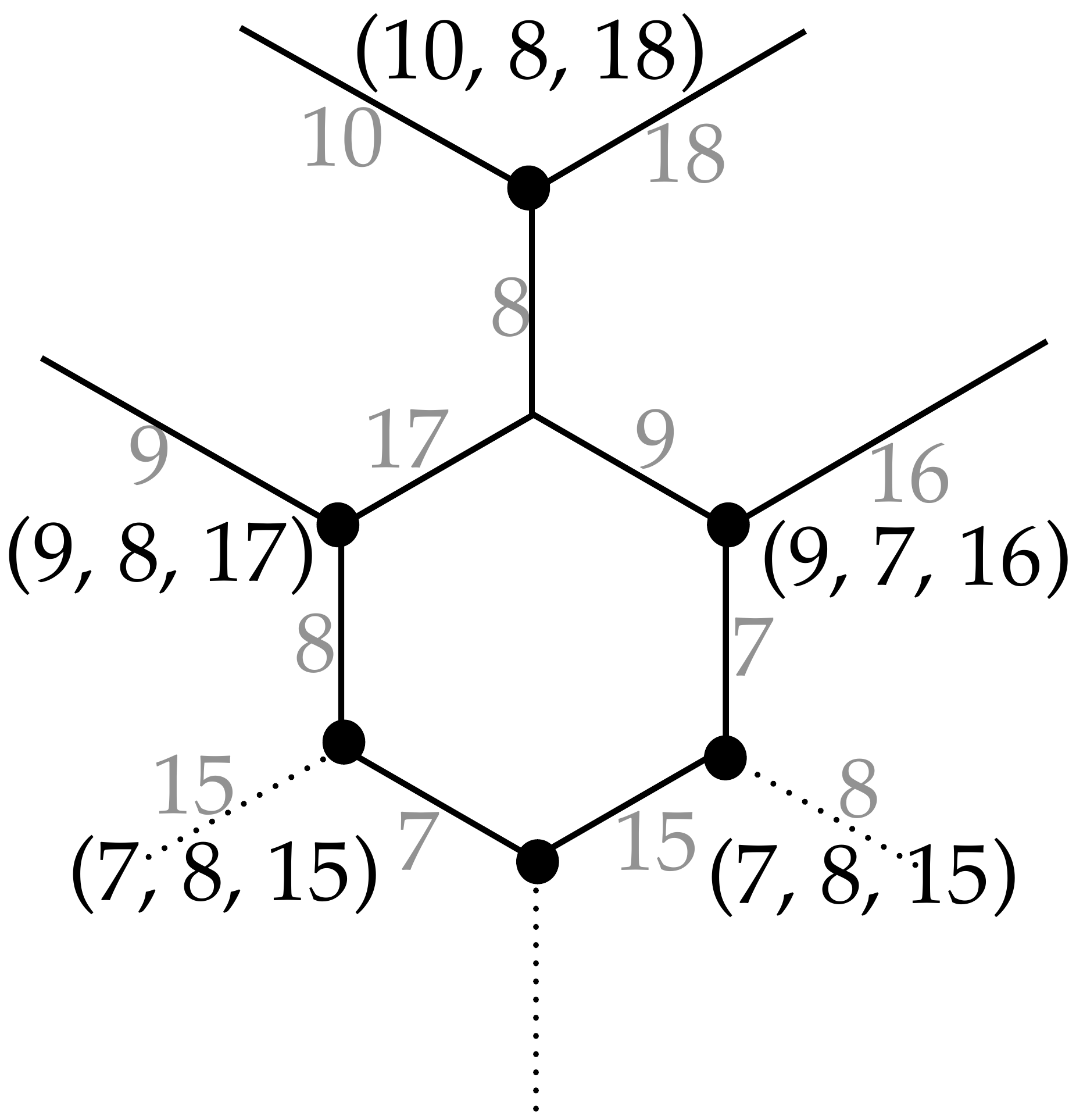}}
\caption{The weighted dual graph viewed in ${\mathbb R}^3$}
\label{Fig18}
\end{figure}

Given the weighted dual graph (coming from the hive, with edge weights given by the positive difference of adjacent hive entries) we can view the edge weights on the edges incident to a given vertex as coordinates on that vertex for a point in a subspace $S$ of ${\mathbb R}^3$, given by $(x,y,z)\in S$ if and only if $x+y=z$.  So, each vertex in the dual graph corresponds to a coordinate in $S$, as shown in Figure~\ref{Fig18} for the top of the weighted dual graph in Figure~\ref{Fig14}.
Note that more than one vertex can correspond to the same point in $S$, as is the case with $(7, 8, 15)$.

This correspondence between vertices of the weighted dual graph and points in a hyperplane in
${\mathbb R}^3$
 creates the {\em honeycomb} for the dual graph. We first plot in $S$ the vertices of the dual graph (possibly with multiplicities) using the coordinates of the edge weights as described above.  Vertices in the honeycomb are connected if the corresponding vertices were adjacent in the dual graph of the hive.  By the flow requirements of vertices of the dual graph, all points $(x,y,z)$ of a honeycomb lie on the hyperplane in ${\mathbb R}^3$ given by $x+y=z$, so we may represent a honeycomb as a planar graph lying in this subspace. The honeycomb for our hive example is shown below.

\begin{figure}[H]
\centering{\includegraphics[scale=.3]{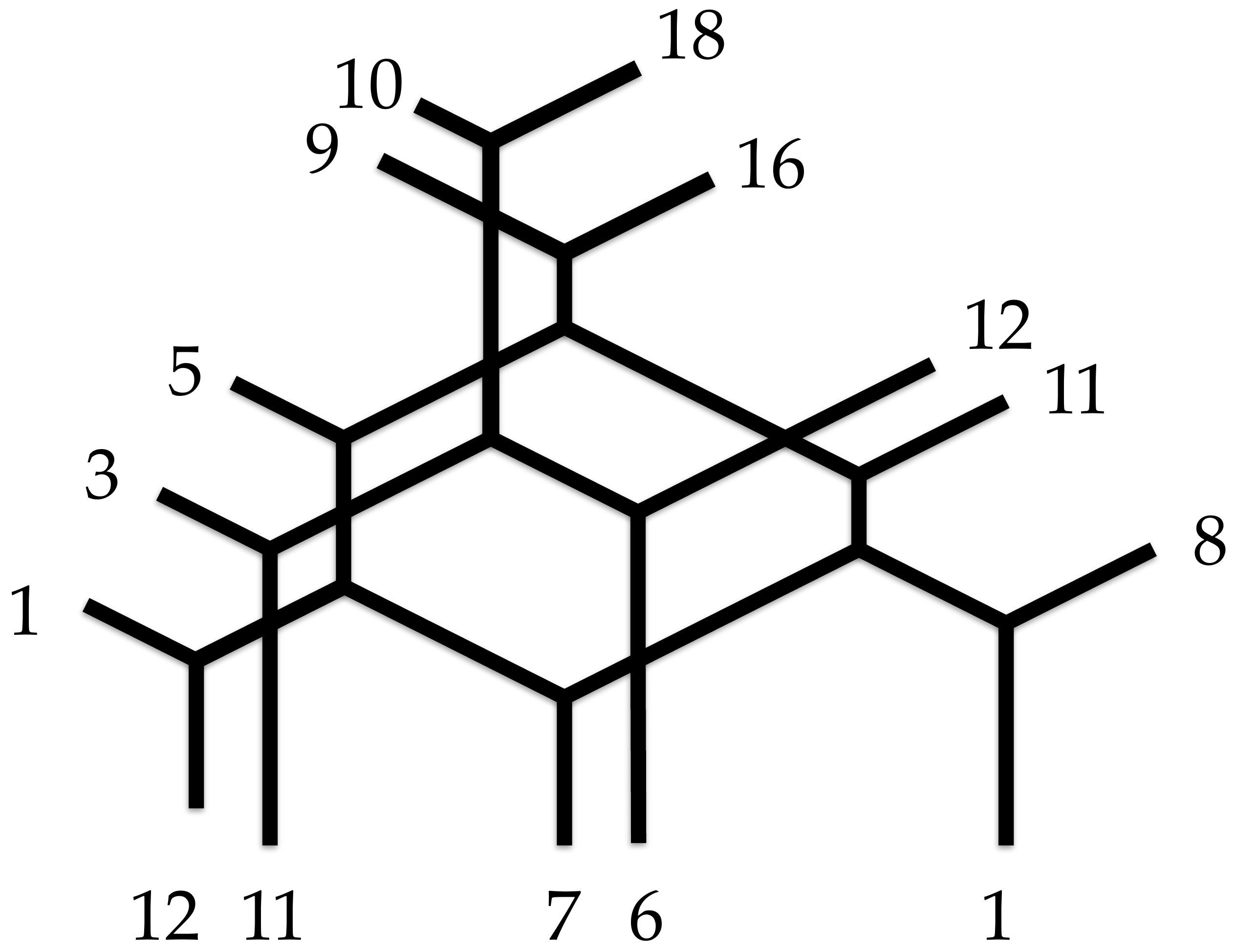}}
\caption{The honeycomb}
\label{Fig20}
\end{figure}

Notice that all the points on a given edge in the honeycomb will have one coordinate that is constant.  For example, the edge labeled 9 starting on the left side of the honeycomb consists of points of the form $(9, y,z)$.  This will allow us to define a capacity for a flow on an edge:  an edge's capacity is defined to be the value of the constant coordinate of the edge.  As with the hive, we define the ``$\mu$ direction" to be along edges going north-west to south-east, the ``$\nu$ direction" to be along edges going south to north and the ``$\lambda$ direction" to be on edges going south-west to north-east.  Note that the capacity of each edge in these three directions is the appropriate part of $\mu$, $\nu$ or $\lambda$, so the ``type" of the honeycomb is the same as the type of the \LR\ filling (and the hive and the weighted dual graph).

Transverse crossings in the honeycomb correspond to adjacent vertices in the dual graph with the same coordinates:
\begin{figure}[H]
\centering{\includegraphics[scale=.3]{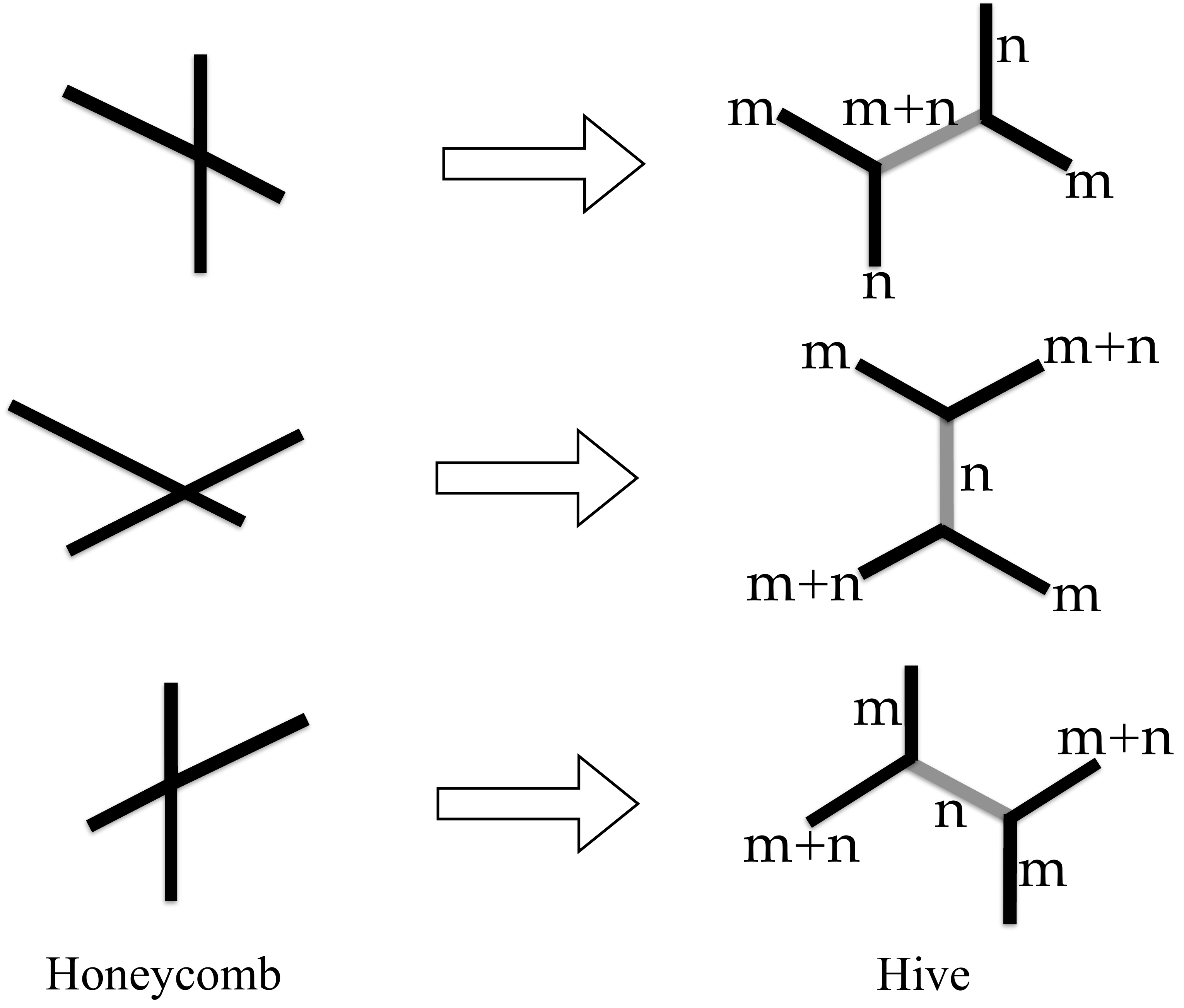}}
\caption{Transverse crossings correspond to dual graph vertices with the same coordinates}
\label{Fig21}
\end{figure}

The canonical flow of the dual graph of a hive translates to a canonical flow on the honeycomb in the obvious manner.  For example, Figure \ref{Fig22} below shows the honeycomb with the $\mu$ flow and the flow of the $1$'s and $2$'s.  This is the same flow as shown on the dual graph in Figure \ref{Fig17}.  The converse holds as well: a canonical flow on a honeycomb produces a \LR\ filling. Knutson and Tao~\cite{knut} show that every honeycomb may be associated to a unique hive, from which a canonical flow on the honeycomb may be determined, and from this flow we may read off the \LR\ filling.

\begin{figure}[H]
\centering{\includegraphics[scale=.3]{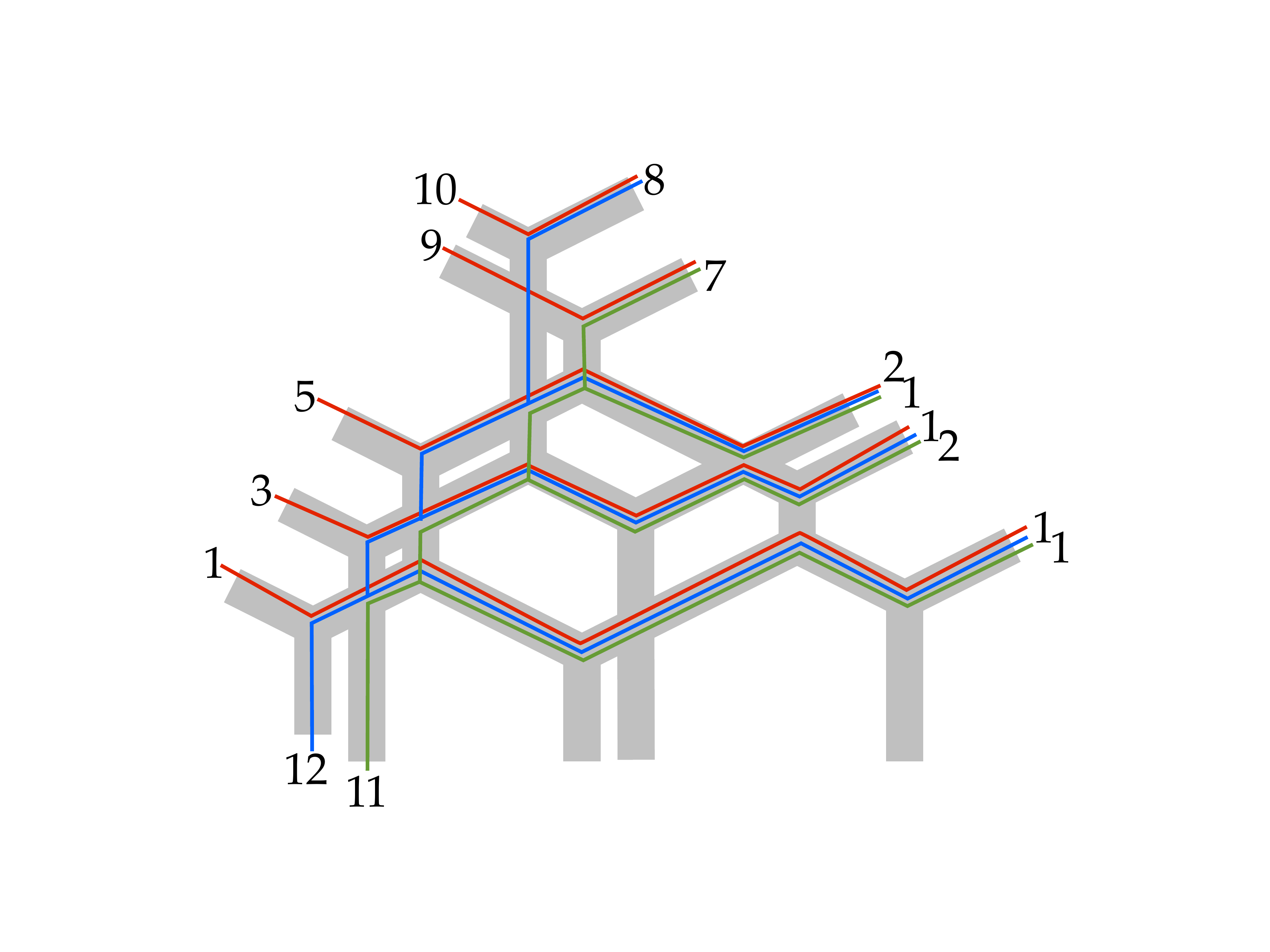}}
\caption{$\mu$ flow and the flow of the $1$'s and $2$'s in the honeycomb}
\label{Fig22}
\end{figure}

\section {Proof that the Algorithm works}

As described above in Section 3, we claimed that our algorithm takes \LR\ fillings of type $(\mu, \nu;\lambda)$ and $(\mu', \nu';\lambda')$, respectively, producing a \LR\ filling of type $(\mu\oplus \mu', \nu\oplus \nu' ;\lambda\oplus \lambda')$.  In this section we prove that the algorithm terminates in a \LR\ filling of the proper type.  We do this by working with the two corresponding flows on the honeycombs associated with the \LR\ fillings of type $(\mu, \nu;\lambda)$ and $(\mu', \nu';\lambda')$.

Our plan is as follows: In~\cite{knut} Knutson and Tao provide a definition of a honeycomb that is independent of the underlying dual graph of a hive, and show that the associated hive may be derived directly from the honeycomb. Furthermore, this definition implies that the overlay of two honeycombs (in the hyperplane $x+y=z$ in ${\mathbb R}^3$) is another honeycomb (indeed, the overlay is one of the ways that multiplicities for edges and vertices may be realized). However, the overlay of the two associated canonical flows is {\em not} typically a canonical flow for the overlay honeycomb.  We shall show that the process of resolving non-canonical flows for the overlay of two honeycombs will match, step by step, the algorithm we presented for sums of \LR\ tableaux, thus proving the algorithm terminates in a \LR\ filling for the sum of the partitions.

So we begin with the two \LR\ fillings.  These correspond to two honeycombs, each with their canonical flow.
\begin{figure}[H]
\centering{\includegraphics[scale=.4]{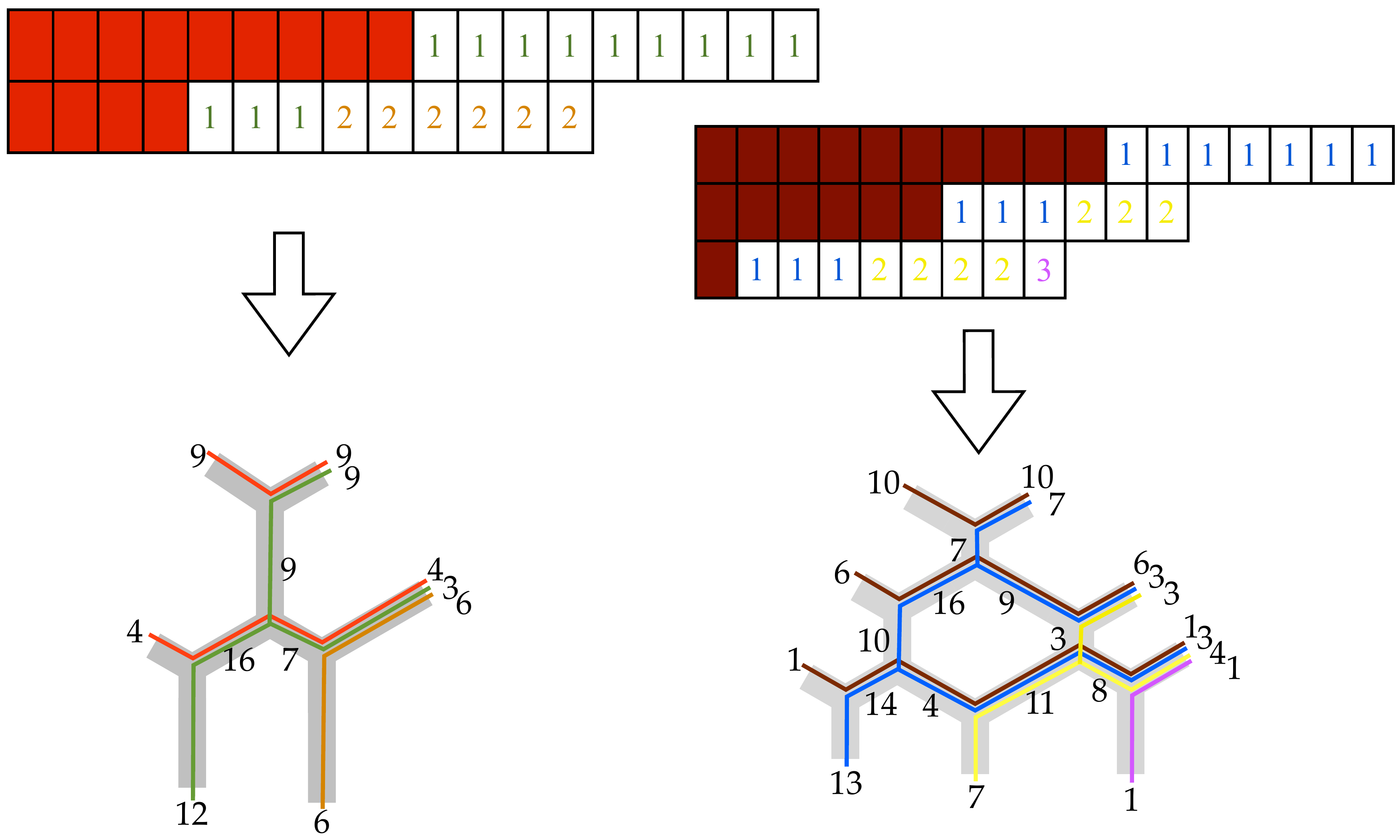}}
\caption{ The two honeycombs associated with the two \LR\ fillings}
\label{Fig24}
\end{figure}
We create the honeycomb that will correspond to the summed tableaux by overlaying the two honeycombs. The associated flows for each honeycomb are also overlayed, producing the ``canonically wrong" flow, which we shall resolve.
\begin{figure}[H]
\centering{\includegraphics[scale=.3]{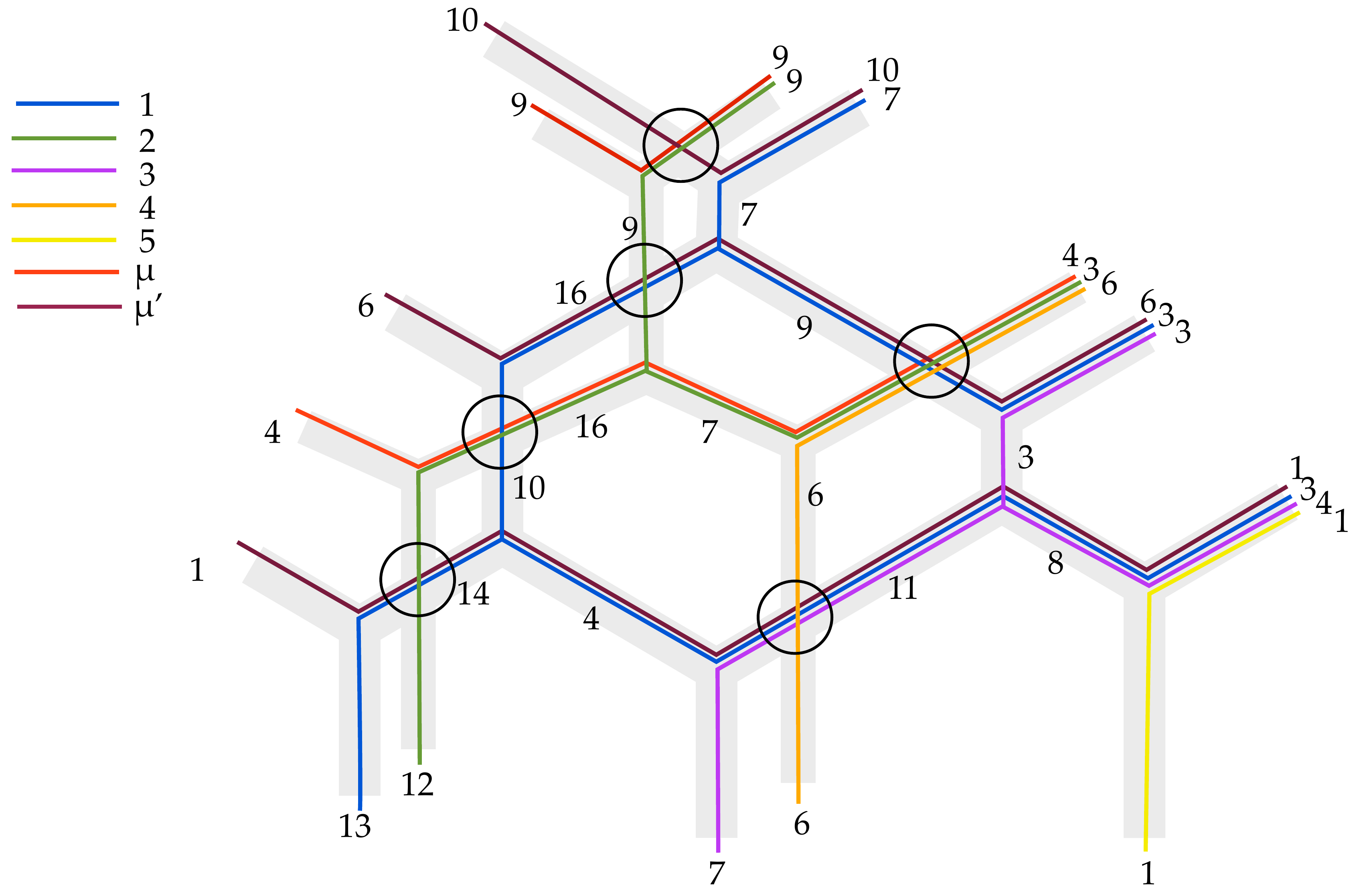}}
\caption{Overlay of two honeycombs and the ``canonically wrong" flow}
\label{Fig26}
\end{figure}

This overlay honeycomb has type $(\mu\oplus \mu', \nu\oplus \nu' ;\lambda\oplus \lambda')$; the content comes from relabeling the content from each diagram, (as was done in the summed diagram) so that the number of $1$'s is greater than or equal to the number of $2$'s, etc.  If two parts of $ \nu\oplus \nu'$ are equal, assign the smaller number to the part that ends in the higher row of $\lambda\oplus \lambda'$.

The overlayed flow will not necessarily be canonical.  In fact, it will have non-canonical flow that mimics exactly the errors in the initial summed filling. (See Figure \ref{Fig2} for an example.  The circled vertices in Figure \ref{Fig26} show the corresponding errors in the honeycomb flow.)  The overlay of honeycombs and their flows can create two types of non-canonical flow:
\begin{figure}[H]
\centering{\includegraphics[scale=.3]{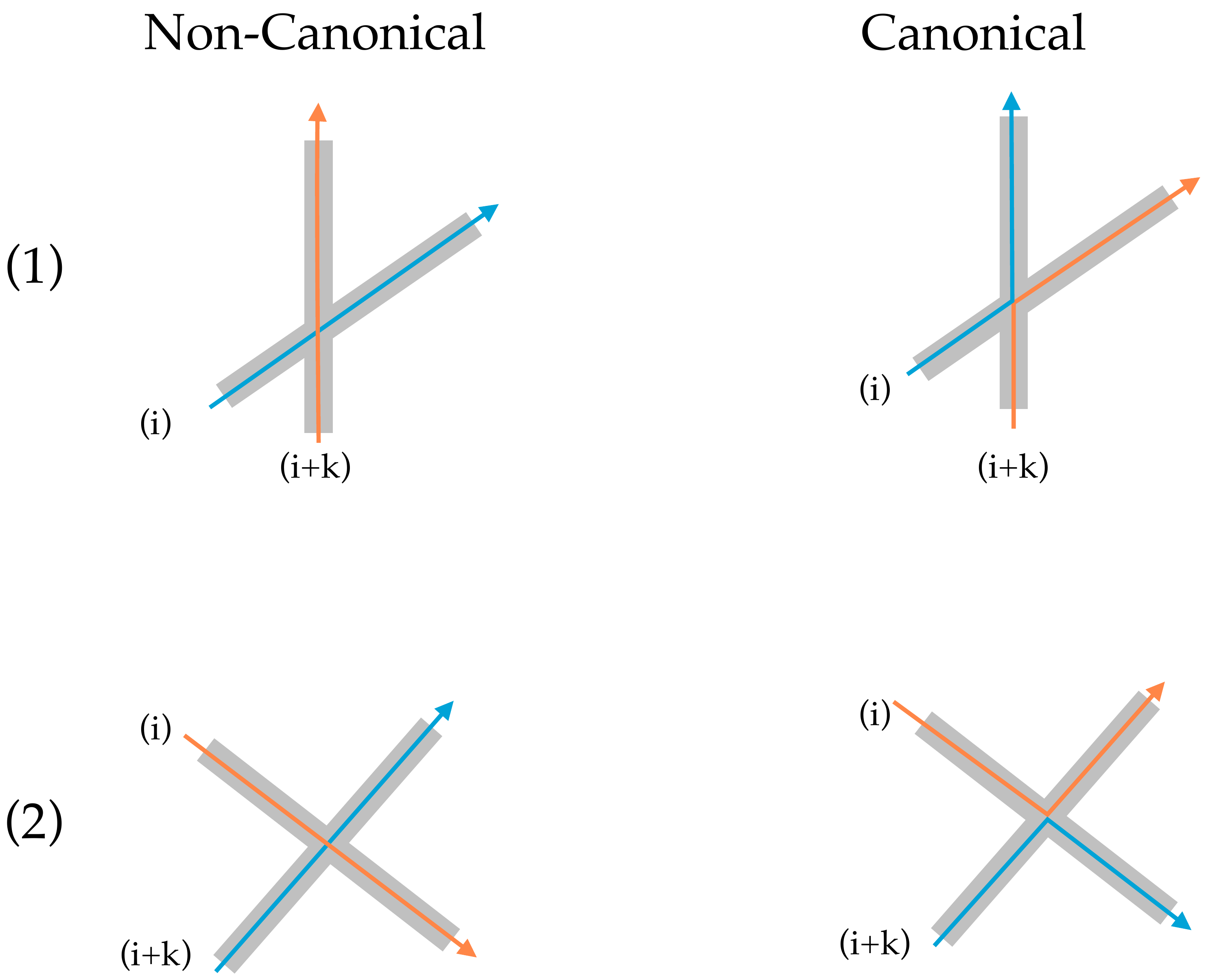}}
\caption{Types of non-canonical flows in the honeycomb overlay}
\label{Fig27}
\end{figure}
We now show how to correct the overlay flow to produce its canonical flow and that this canonical flow is indeed the one associated to the corrected \LR\ filling of the sum of the two diagrams. This will prove that the algorithm given in Section 3 terminates in a \LR\ filling of the summed diagram.

First, the two types of non-canonical flows that can be seen in the overlay flow (Figure \ref {Fig27}) correspond to the following types of errors in the sum filling:  Type (1) flow errors are ``$i$ in row $i$" or word violations; Type (2) flow errors are ``$\mu$ switching" or column-strict violations.  There are no $\lambda$ violations in the overlay; the parts of $\lambda\oplus\lambda'$ are ordered from largest to smallest in the honeycomb. Although it is unimportant to the final, corrected flow, to be consistent with the order in which errors are corrected in the \LR\ diagram, we will correct each type of non-canonical flow from top to bottom, east to west, as encountered.

A ``$\mu$ switching" error occurs when the honeycomb overlay produces a transverse crossing  of the $\mu$ flow, as shown:
\begin{figure}[H]
\centering{\includegraphics[scale=.3]{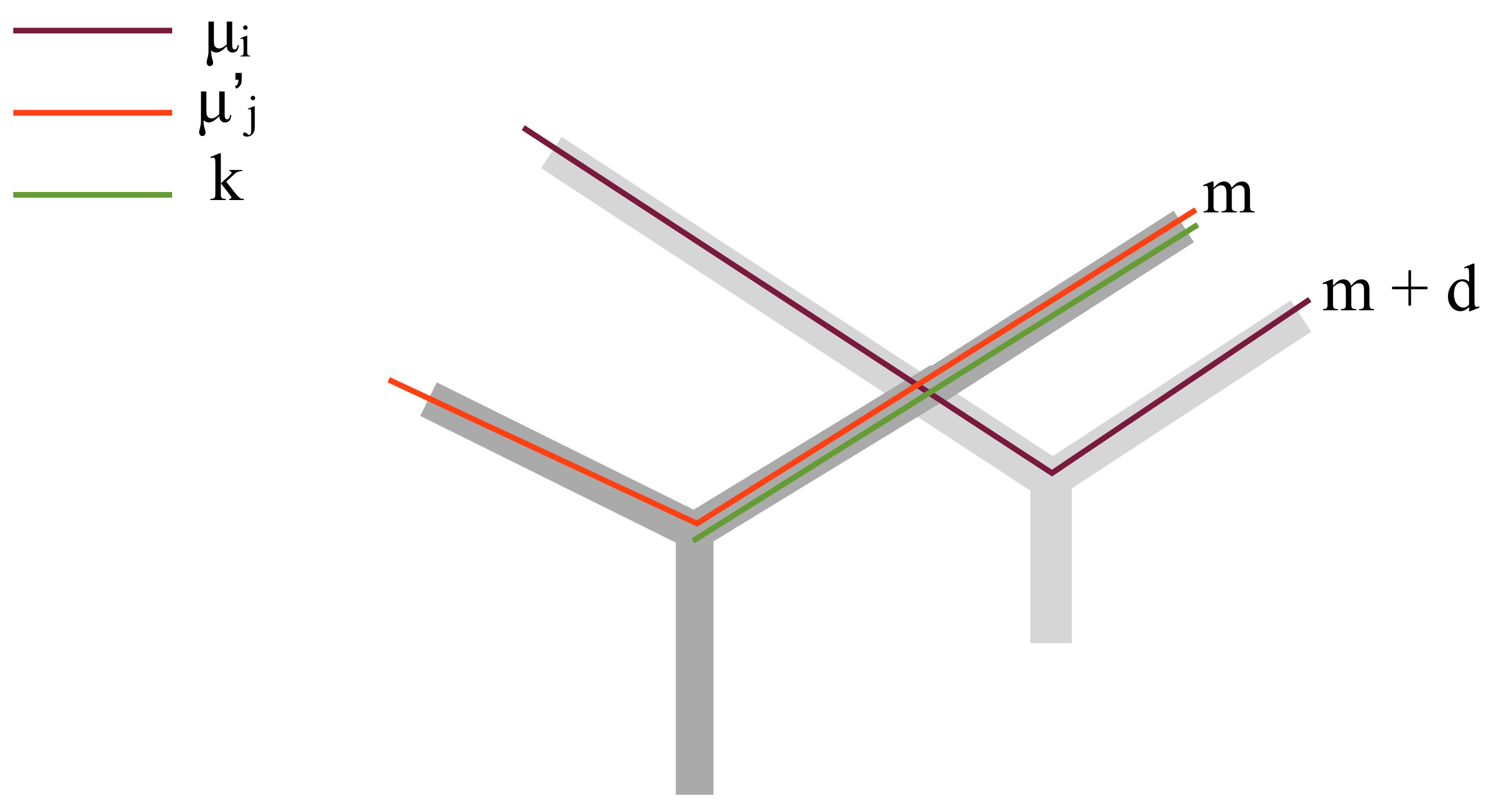}}
\caption{$\mu$ switching error in flow}
\label{Fig29}
\end{figure}

The two parts of $\mu \oplus \mu'$ come from different honeycombs, with the longer part crossing below the smaller part.  We correct this error by mimicking the process done on the \LR\ diagrams.  Let $\mu_i = m+d$ and $\mu'_j = m$.  We need to swap $m+d$ units of flow on the upper edge, which will consist of $\mu$ flow and $\nu$ flow, with $m+d$ units of flow on the lower edge, which will be entirely $\mu$ flow. Note that this is always possible because each edge's label (and hence the total flow on that edge) corresponds to its $z$-coordinate, which is always at least as big as its $x$ coordinate, the $\mu$ flow value.

So we correct the crossing by taking $m+d$ units of  flow ($\mu$, 1's, 2's, etc. in order) on the upper edge and swapping that with $m + d$ units of flow on the lower edge. For example, we swap the following.  Note that this agrees completely with the ``$\mu$ flow" correction on the sum of \LR\ diagrams.

\begin{figure}[H]
\centering{\includegraphics[scale=.3]{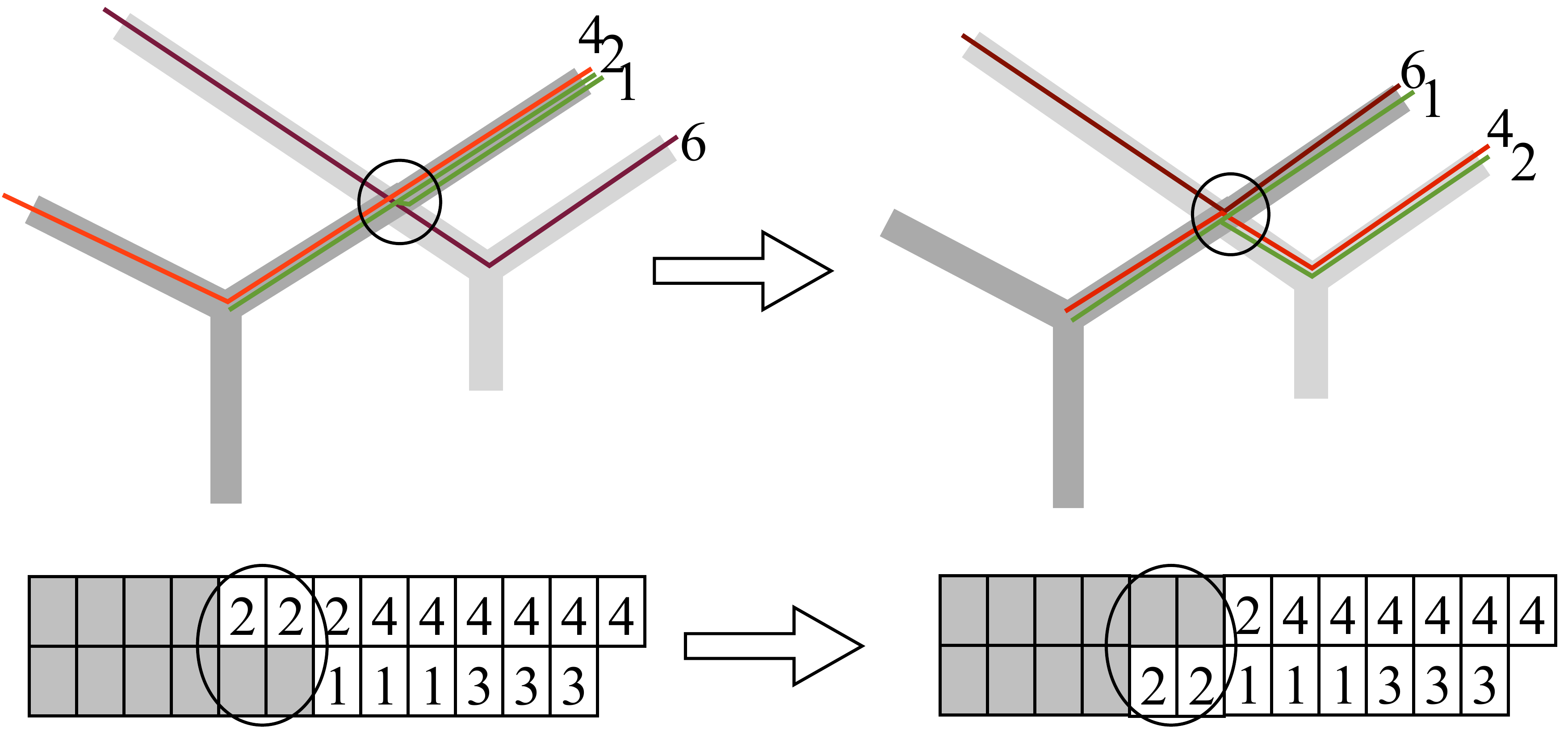}}
\caption{Strand swapping to correct the $\mu$ flow error}
\label{Fig30}
\end{figure}

After all $\mu$ switching errors are fixed (our example has one more at the top circled crossing in Figure \ref {Fig26}), we correct ``$i$ in row $i$" errors.  These errors occur when, in overlaying the two honeycombs, an $i+k$ ends in row $i$. In the figure below, $k=1$.
\begin{figure}[H]
\centering{\includegraphics[scale=.3]{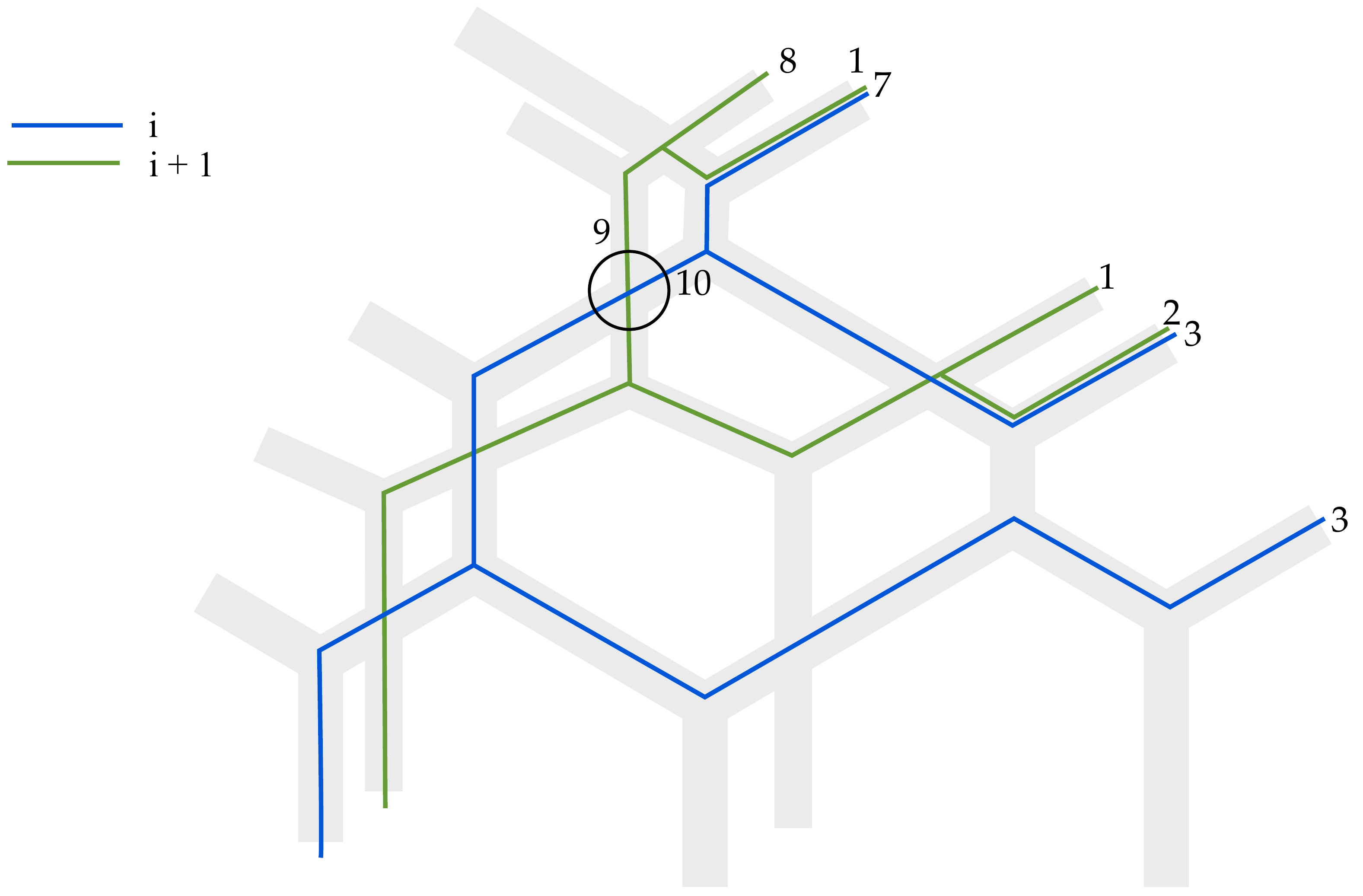}}
\caption{``$i$ in row $i$ error in the overlay of the honeycombs}
\label{Fig31}
\end{figure}
We correct this error much as in the previous case, by swapping equivalent amount of flow at the ``bad" intersection, so that the flow is canonical there.  First as shown below, this sort of transverse crossing comes from two duplicate vertices in the dual graph. (Note that the green flow on the edge labeled ``$\ell -n$" in the dual graph consists of $n$ units flowing {\em backwards}, which is how non-canonical flows arise on the dual graph, and one reason analyzing flows on honeycombs is simpler.)
\begin{figure}[H]
\centering{\includegraphics[scale=.4]{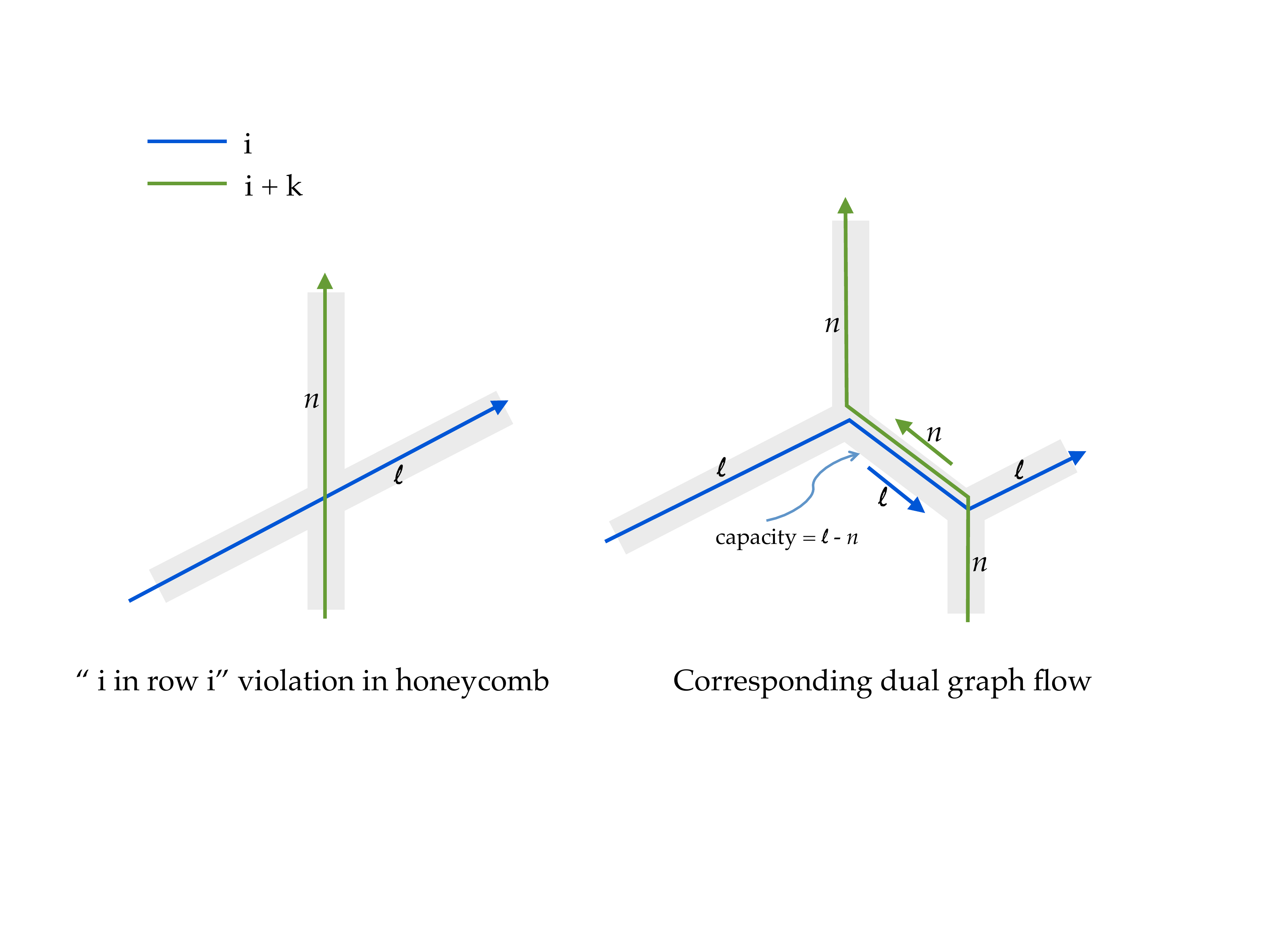}}
\caption{Dual graph flow}
\label{Fig32}
\end{figure}
Recalling that the sum of the flows in the $\mu$ and $\nu$ directions equals the flow in the $\lambda$ direction, we see that $\ell \ge n$. So we exchange $n$ units of $i$ flow in the $\lambda$ direction with $n$ units of $i+k$ flow in the $\nu$ direction.  For example, to fix the error at the circled vertex, we need to swap $n = 9$ units of $2$-flow in row 1 with $1$-flow in row 2.  The swapped flow is depicted in a heavier line.
\begin{figure}[H]
\centering{\includegraphics[scale=.3]{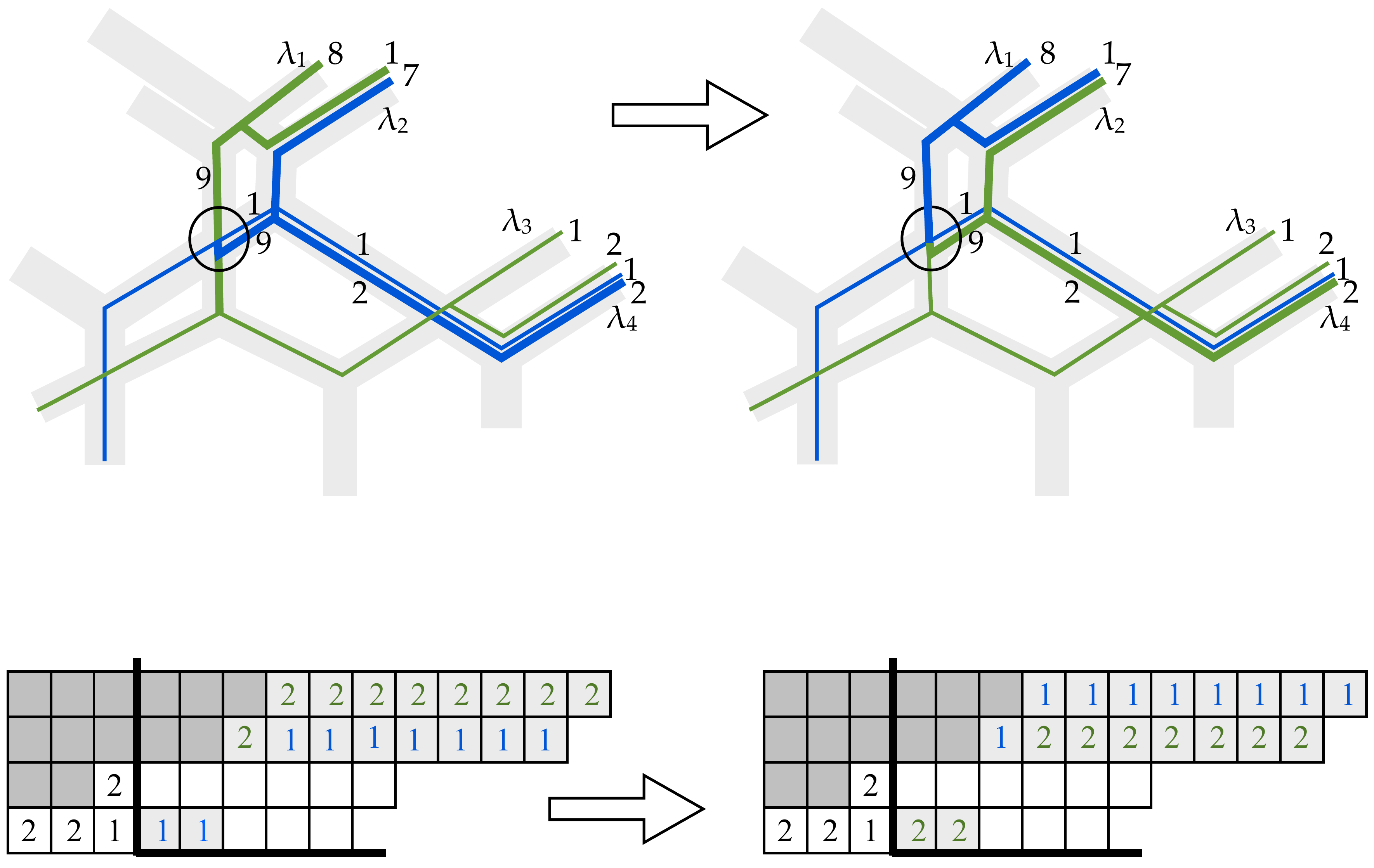}}
\caption{The flow, before and after the swap}
\label{Fig33}
\end{figure}
Note that this accomplishes the swap done on the \LR\ diagram to fix the same problem, as seen above.  Indeed, as we can see in Figure \ref{Fig33}, this type of transverse crossing forces at least as many $i$'s to flow north-east (in the $\lambda$ direction) as there are $(i+1)$'s flowing north (in the $\nu$ direction), corresponding to the first row of the \LR\ diagram in which the strand of $i$'s is at least as long as the strand of $(i+1)$'s.

\begin{figure}[H]
\centering{\includegraphics[scale=.3]{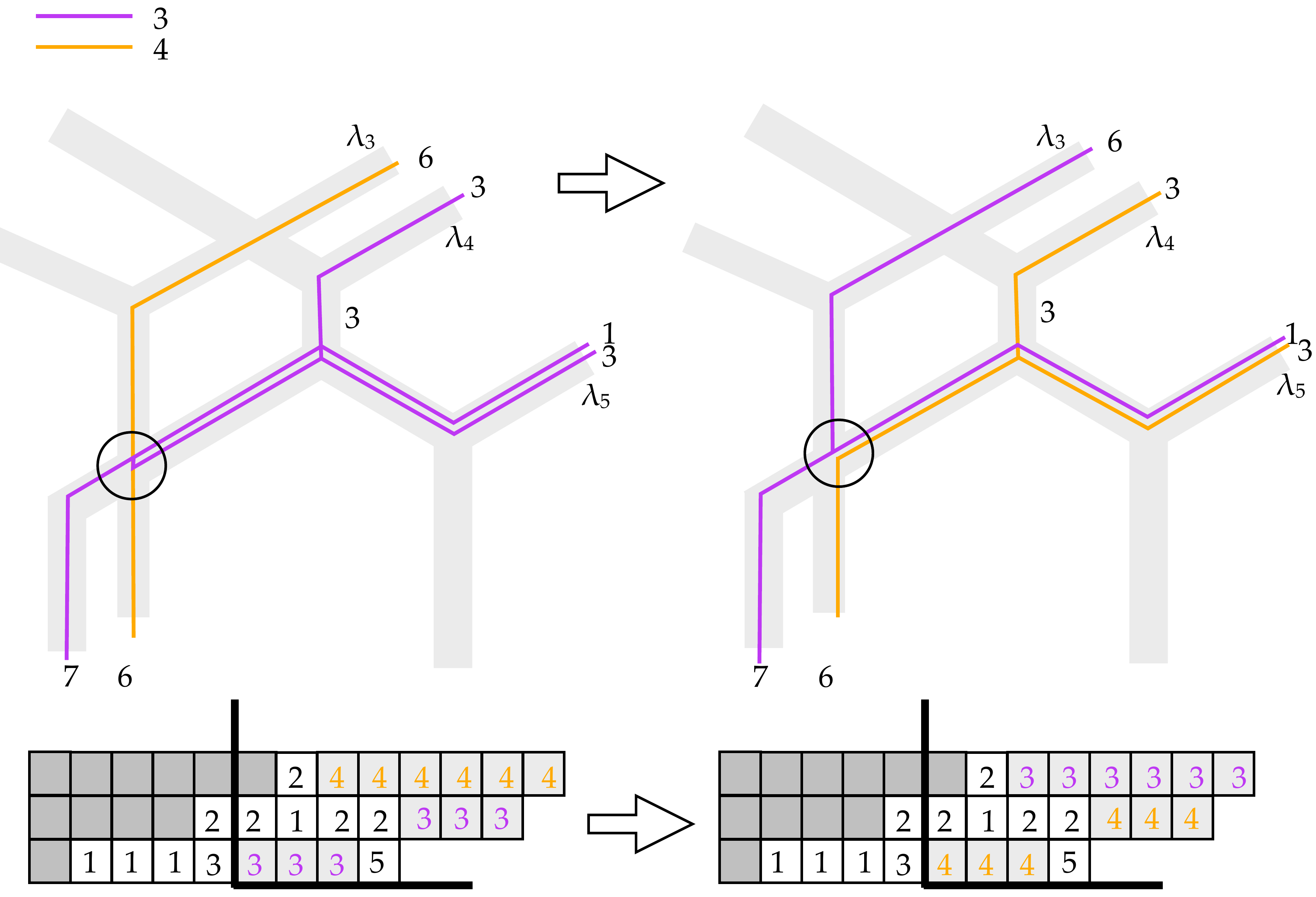}}
\caption{``3 in row 3" flow fixed, and the corresponding rows of the \LR\ diagram}
\label{Fig34}
\end{figure}

Below, Figure \ref {Fig34} shows the flow swap on the ``3 in row 3" problem, indicated by the non-canonical flow circled.

Continuing to follow the algorithm described in Section 3, we now tackle word or column-strict flow violations.  Word violations correspond to Type (1) non-canonical flows (see Figure \ref {Fig27}).  Column-strict violations correspond to Type (2).  As such, these errors are corrected in exactly the same manner as a $\mu$-flow (another Type (1) error) and ``$i$ in row $i$" (a Type (2) error).

\begin{figure}[H]
\centering{\includegraphics[scale=.25]{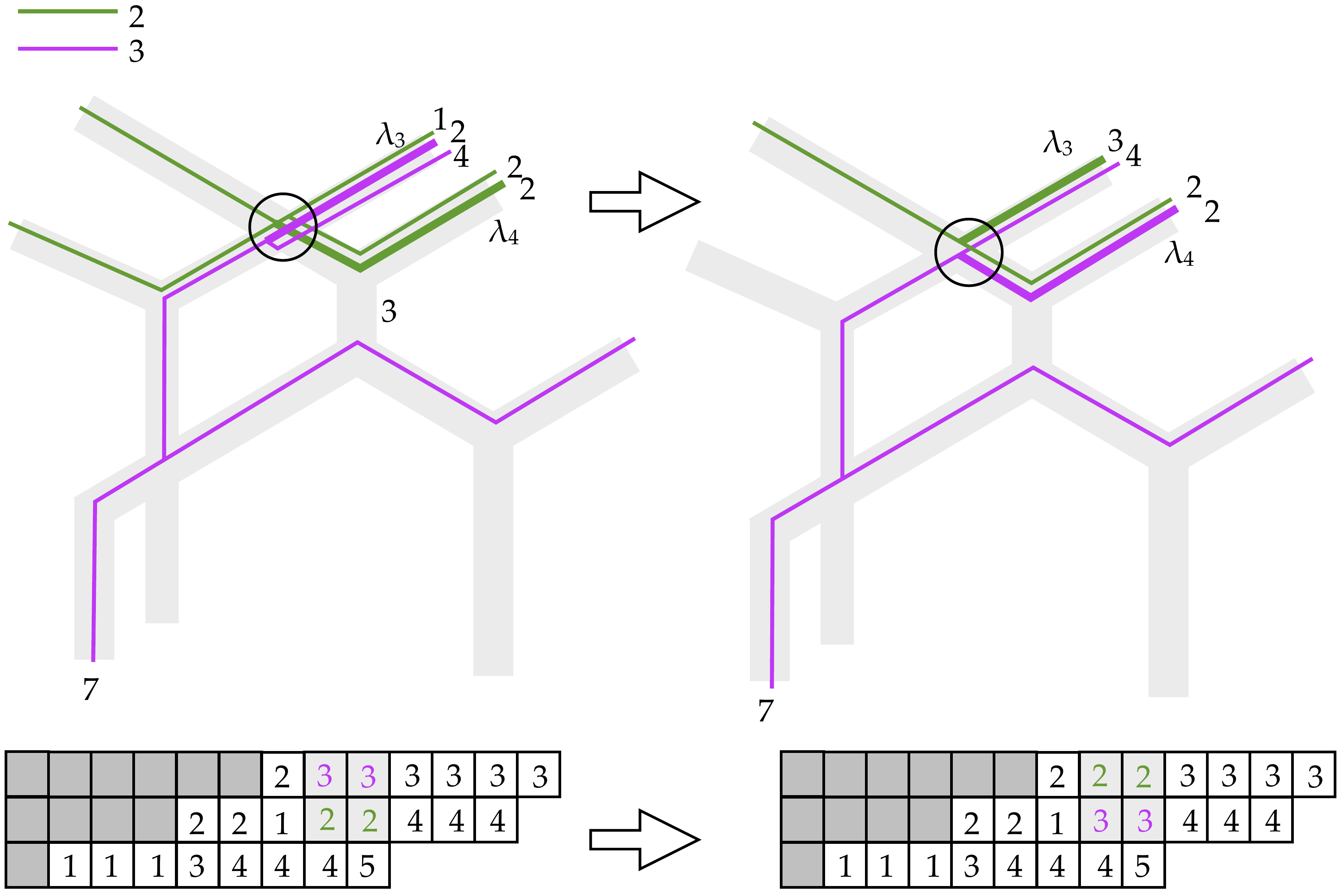}}
\caption{Column-strict flow correction, and the corresponding rows of the \LR\ diagram}
\label{Fig35}
\end{figure}

For example, Figures \ref {Fig6} and \ref {Fig7} show the \LR\ diagram with a column-strict violation fixed.  The corresponding non-canonical flow and strand swap are shown in Figure \ref{Fig35}.
After all non-canonical flow intersections are corrected, as described above, the resulting flow is canonical, and hence, the corresponding diagram, the one produced by the algorithm in Section 3, is indeed \LR, save for any row strict violations that may remain in the diagram. However, we may, as claimed, resolve those violations without risk of producing any new errors in the diagram since the counts of $1$'s, $2$'s, etc., in the diagram will remain unchanged, and will match the counts of the necessarily canonical flow on the right edge of the honeycomb. Hence, re-ordering within the row can only result in a filling matching that of the canonical flow.

Finishing our example, the corrected, canonical, flow on the honeycomb overlay is:
\begin{figure}[H]
\centering{\includegraphics[scale=.3]{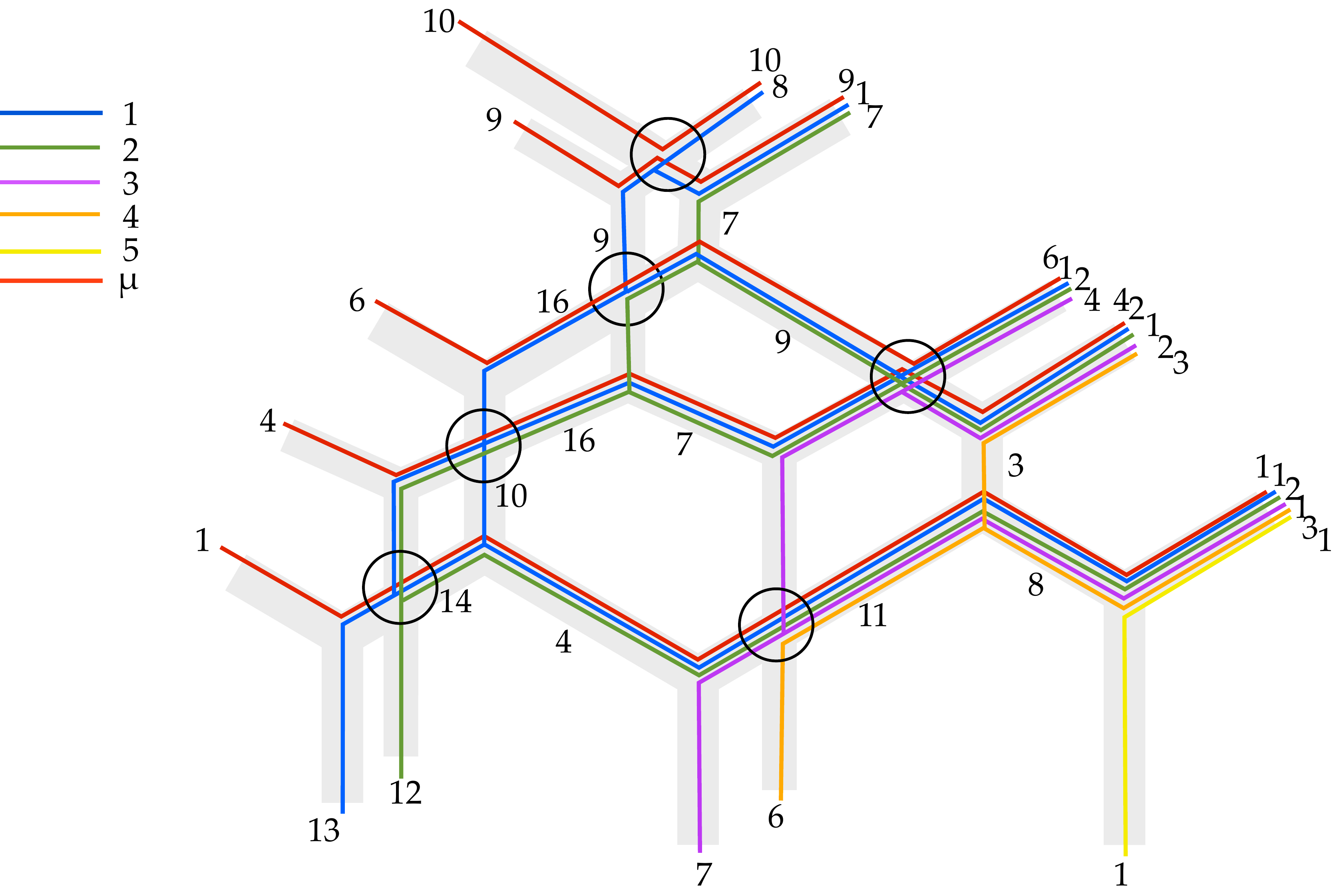}}
\caption{The corrected, canonical flow of the overlay}
\label{Fig36}
\end{figure}
This flow corresponds to the corrected \LR\ diagram in Figure \ref{Fig9}.
\section {Conclusion}
As mentioned earlier, the work of King, Tollu, and Toumazet~\cite{KTT} gave conditions that determined when a \LR\ coefficient may be factored as a product of two other \LR\ coefficients. In our terminology, this would amount to finding conditions such that
\[ c_{\mu \oplus \mu', \nu \oplus \nu'}^{\lambda \oplus \lambda'} = c_{\mu \nu}^{\lambda} \cdot c_{\mu' \nu'}^{\lambda'}. \]
The results of~\cite{KTT} demonstrated how one might decompose a \LR\ filling (indexed by $c_{\mu \oplus \mu', \nu \oplus \nu'}^{\lambda \oplus \lambda'}$) and show how it was the result of a specialized ``overlay" of hives associated to the filling. Our results show, conversely, how to recover any filling in $LR(\mu \oplus \mu', \nu \oplus \nu' ; \lambda \oplus \lambda')$ (determined by an essential triple or not) as a sum of fillings, one from $LR(\mu,\nu; \lambda)$ and one from $LR(\mu', \nu' ; \lambda')$. As the honeycomb example at the start of our paper showed, there are honeycomb decompositions that do not satisfy the conditions of~\cite{KTT}. Indeed, by our algorithm for the sum of \LR\ fillings, there is always map:
\[ LR(\mu,\nu; \lambda) \times LR(\mu', \nu' ; \lambda') \rightarrow LR(\mu \oplus \mu', \nu \oplus \nu' ; \lambda \oplus \lambda'). \]
As described in the introduction, King, Tollu, and Toumazet~\cite{KTT} provide conditions under which the above map is a bijection, but often it fails to be onto, and sometimes it is not one-to-one. Such examples suggest that this phenomenon may be interpreted as a necessary property for points (\LR\ fillings) on the vertices of the polytope of \LR\ fillings associated to a given triple. The ability to deform a filling in more than two independent directions in this polytope is reflected in that filling having not only one, but more than one decomposition. Further, our methods, based (ultimately) on the combinatorics of honeycombs, demonstrate that some problems in the geometry of \LR\ fillings may be analyzed more easily using honeycombs, instead of hives as has more often been the case: The ``strand swapping" at the heart of the algorithm is rather difficult to describe on a hive and is still difficult to work with on the dual graph, but the process becomes transparent when represented on the honeycomb.

Another setting in which \LR\ fillings have been applied is in the study of invariant factors of matrices over rings with a discrete, and even a real-valued, valuation. In these cases, a triple of partitions $(\mu, \nu, \lambda)$ is associated to square matrices $M$ and $N$ over a valuation ring $R$ such that the sequence of {\em orders} of the invariant factors (with respect to the valuation) is $\mu$ for the matrix $M$, $\nu$ for the matrix $N$, and $\lambda$ for the product $MN$. It has been shown~\cite{fulton} that the {\em same} triples of partitions (possibly real-valued) appear in this matrix setting (possibly using a ring with real valuation) as those appearing in the Hermitian case. In the valuation ring setting the present authors have determined methods to associate a \LR\ filling to matrix pairs, and conversely~\cite{me,lrcon,lrreal}. Furthermore, calculations by the present authors corroborate that \LR\ fillings associated to the direct sum of matrices in the valuation ring context do, in fact, correspond to the sum of the fillings of the summands, as calculated by our algorithm here. We feel that analyzing the summation algorithm in the ring theoretic context may contribute to our understanding in how to construct \LR\ fillings and/or honeycombs associated to matrix pairs in the Hermitian context.


\end{document}